\documentclass{article}

\usepackage[utf8]{inputenc}

\usepackage[english]{babel}

\usepackage[letterpaper,top=2cm,bottom=2cm,left=3cm,right=3cm,marginparwidth=1.75cm]{geometry}
\usepackage{siunitx} 
\usepackage{amsthm,amsmath,amsfonts,amssymb, mathtools}
\usepackage{graphicx}
\usepackage[colorlinks=true, allcolors=blue]{hyperref}
\usepackage{xcolor}
\usepackage[anythingbreaks]{breakurl}
\usepackage{soul}

\usepackage[ruled,vlined]{algorithm2e}
\usepackage{lipsum}
\usepackage{authblk} 

\usepackage{url}

\usepackage{natbib}

\usepackage{float}
\usepackage{orcidlink}
\pdfoutput=1

\usepackage{xr}
\mathtoolsset{showonlyrefs}

\title{Adaptive estimation for Weakly Dependent Functional Times Series}
\author[1,2]{Hassan Maissoro}
\author[1]{Valentin Patilea}
\author[1]{Myriam Vimond\thanks{Corresponding author. E-mail: myriam.vimond@ensai.fr} \orcidlink{0000-0002-7722-9933}}

\date{January 29, 2024}

\affil[1]{\small Univ Rennes, Ensai, CNRS, CREST (Center for Research in Economics and Statistics) - UMR 9194, F-35000 Rennes, France}
\affil[2]{\small Datastorm, F-91120 Palaiseau - France}

\newcommand{\Cc}{\mathcal{C}}
\newcommand{\EE}{\mathbb{E}}

\newcommand{\Hh}{\mathcal{H}}

\newcommand{\Ll}{\mathcal{L}}
\newcommand{\LL}{\mathbb{L}}
\newcommand{\NN}{\mathbb{N}}
\newcommand{\Oo}{\mathcal{O}}
\newcommand{\PP}{\mathbb{P}}

\newcommand{\RR}{\mathbb{R}}

\newcommand{\Xx}{\mathcal{X}}

\newcommand{\ZZ}{\mathbb{Z}}

\newcommand{\vertiii}[1]{{\left\vert\kern-0.25ex\left\vert\kern-0.25ex\left\vert #1 
    \right\vert\kern-0.25ex\right\vert\kern-0.25ex\right\vert}}

\RequirePackage{enumerate}
\RequirePackage{comment}
\RequirePackage{ifthen}
\RequirePackage{subcaption}
\RequirePackage{siunitx}
\RequirePackage[ruled,vlined]{algorithm2e}
\usepackage{dsfont}

\theoremstyle{plain}
\newtheorem{theorem}{Theorem}

\newtheorem{definition}{Definition}
\newtheorem{lemma}{Lemma}

\newtheorem{example}{Example}
\newtheorem{remark}{Remark}

\DeclareSIUnit{\octet}{o}
\newcounter{assumptionHt}
\newcounter{assumptionH}
\newcounter{assumptionE}
\newcounter{assumptionD}
{%
\begin{enumerate}[({G}1)]%
    \setcounter{enumi}{\value{assumptionHt}}%
    }{%
    \setcounter{assumptionHt}{\value{enumi}}%
\end{enumerate}
}
\newenvironment{assumptionH}%
{%
\begin{enumerate}[({H}1)]%
    \setcounter{enumi}{\value{assumptionH}}%
    }{%
    \setcounter{assumptionH}{\value{enumi}}%
\end{enumerate}
}
{%
\begin{enumerate}[({E}1)]%
    \setcounter{enumi}{\value{assumptionE}}%
    }{%
    \setcounter{assumptionE}{\value{enumi}}%
\end{enumerate}
}
{%
\begin{enumerate}[({D}1)]%
    \setcounter{enumi}{\value{assumptionD}}%
    }{%
    \setcounter{assumptionD}{\value{enumi}}%
\end{enumerate}
}
\newcommand{\assrefHt}[1]{(\hyperref[#1]{G\ref{#1}})}
\newcommand{\assrefH}[1]{(\hyperref[#1]{H\ref{#1}})}
\newcommand{\assrefE}[1]{(\hyperref[#1]{E\ref{#1}})}
\newcommand{\assrefD}[1]{(\hyperref[#1]{D\ref{#1}})}

\newcommand{\EEn}{\mathbb{E}_{n}}
\newcommand{\EEMT}{\mathbb{E}_{M,T}}

\newcommand{\n}[1]{{n_{#1}}}

\newcommand{\Xtemp}[1]{%
	\ifthenelse{\equal{#1}{0}}
		{X^{(\n0)}}
		{
			\ifthenelse{\equal{#1}{1}}
				{X^{[\n1]}}
				{X^{(#1)}}
		}
	}
\newcommand{\hXtemp}[1]{%
	\ifthenelse{\equal{#1}{0}}
		{\widehat X^{(\n0)}}
		{
			\ifthenelse{\equal{#1}{1}}
				{\widehat X^{[\n1]}}
				{\widehat X^{(#1)}}
		}
	}

\newcommand{\Yni}{{Y_{n,i}}}

\newcommand{\Tni}{{T_{n,i}}}

\newcommand{\e}{{\varepsilon}}

\DeclareMathOperator*{\argmin}{arg\,min}

\begin{document}

\maketitle

\begin{abstract}
We study the local regularity of weakly dependent functional time series, under   $L^p-m-$appro\-ximability assumptions. The sample paths are observed with error at possibly random, design points. 
Non-asymptotic concentration bounds of the regularity estimators are derived. 
As an application, we build nonparametric mean and autocovariance functions estimators that adapt to the regularity of the sample paths and the design which can be sparse or dense. We also derive the asymptotic normality of the adaptive mean function estimator which allows for honest inference for irregular mean functions. An extensive simulation study and a real data application illustrate the good performance of the new estimators. 
\end{abstract}



%
%
%
%

%
%
%
%
%
%
%

\section{Introduction}
Functional Data Analysis (FDA) refers to the case where the  observation units are the whole curves (also called trajectories or sample paths). The data set then consists of a collection of $N$ trajectories, modeled by a same  stochastic process defined over some domain.
Dependent functional data arise in fields such as environment \citep{aue2015prediction}, energy \citep{chen2021review}, biology \citep{stoehr2021detecting} or clinical research \citep{martinez2021nonparametric,li2023statistical}. They are often collected sequentially at regular time intervals (e.g. days, weeks) and exhibit a serial dependence. Functional time series (FTS) analysis aims to understand the serial dependence between curves and their dynamics over time.
Several types of dependence for functional data have been studied, such as cumulant mixing conditions, strong mixing, physical dependence, $\LL^p-m-$approximability. See, for example, \citet{hormann2012functional,panaretos2013fourier,chen2015simultaneous,rubin2020sparsely} and their references. 
We consider  FTS that are $\LL^p-m-$approximable, \emph{i.e.}, satisfying a general moment-based notion of weak dependence involving $m-$dependence \citep[see][]{hormann2012functional}.

Most of the textbooks and many FDA articles consider the sample paths observed without error at each point in the domain. In this case, the FDA permits straightforward nonparametric approaches (such as the empirical mean and covariance function estimators), for which an elegant theory is derived based on limit theorems for Hilbert space variables. See, for example, \citet{horvath2012inference}. In real data problems, the curves are only observed by a finite number of noisy measurements, at observation design (or domain) points that are not necessarily regular or identical from one curve to another. Two cases are usually studied: the points of the domain where the curves are observed are the same for all the curves (common design), or they are completely different (independent design). The two situations are different in nature and usually lead to different theoretical results. With a common design there is no information about the stochastic process  between the  design points. 

A common practice in FDA is to first create smoothed curves (usually called functional data objects) from the data points on each curve separately, and then proceed as when the sample paths were observed without error everywhere in the domain. For each curve separately, smoothed curves can be constructed by nonparametric smoothing (splines, kernel smoothing, etc.), or simple linear interpolation. There is no reason however why constructing smoothed curves ignoring the other curves generated by the same stochastic process, should always be an appropriate way to proceed with the FDA. Alternatively, for example, for the mean and covariance function estimation, one can pool all the data points and proceed with nonparametric procedures. See \cite{zhang2016sparse} for independent functional data and \cite{rubin2020sparsely} for the dependent sample paths. However, while pooling the data points of all the curves appears to be effective for independent curves, in a time series context it removes the information about the stochastic dependence between the sample paths. 

Nonparametric methods with  separately smoothed curves, sometimes called `smooth first, then estimate' approaches, are usually recommended for curves observed over a \emph{dense} set of domain points, while pooling is preferred for \emph{sparse} functional data \citep[see][]{yao2005functional,zhang2016sparse}. It is worth noting that the definition of sparse and dense regimes depend on the regularity of the sample paths \citep[see][]{zhang2016sparse}. Furthermore, the minimax convergence rates for the nonparametric methods are expected to depend on the regularity of the sample paths \citep[see][for the mean function estimation]{cai2011}. While the regularity of the sample paths, an intrinsic property of the process generating the functional data, has a major impact on nonparametric methods in FDA, it appears that  little effort has been devoted so far to estimating this regularity and to constructing adaptive methods. In most cases, the sample paths are assumed to have a certain regularity, e.g. twice continuously differentiable. However, many applications produce irregular curves, such as photovoltaic or wind power generation which depend on natural phenomena. 

In this paper we first study an estimation method for the local regularity of the process generating the FTS. In the case of non-differentiable sample paths, the local regularity is given by the local Hölder exponent and constant of the mean squared increments of the process. In the case of differentiable sample paths, the increments of the largest order derivative of the sample paths are considered instead. Our contribution extends that of \citet{golovkine2022learning} who have studied independent and identically distributed functional data. The definition of local regularity considered below is closely related to the notion of local intrinsic stationarity introduced by \cite[page 2060]{hsing2016}. Similar concepts of local regularity are common in continuous-time processes, \citep[see, for example,][Section 5, and their references]{Jirak17}, where the regularity estimation is usually based on a single sample path. 
Using the local regularity estimators and kernel smoothing of the curves, we propose adaptive mean and autocovariance function estimators. The local bandwidths are data-driven, chosen by minimizing explicit quadratic risk bounds. 

The paper is organized as follows. Section~\ref{sec:FTS_model} presents the statistical model associated with the observation of the FTS at discrete points in the domain, in the presence of additive heteroscedastic noise. The local regularity assumptions and the type of weak dependence assumption considered are next introduced. Section~\ref{sec:local-reg-estimator} presents the estimators and their concentration bounds. Both non-differentiable and differentiable sample paths cases are considered. Our estimator is related to the estimation of the Hurst function of a multifractional Brownian motion, and to other regularity parameters studied in the stochastic process theory. See for example \cite{Hsing20}. As an application of our local regularity estimators, in Section \ref{sec:adapt_mean_acov} we propose adaptive kernel estimators  for the mean function and the autocovariance functions, for which we derive the pointwise convergence rates. They adapt to the regularity of the sample paths and the design which can be independent or common, sparse or dense. We also prove the asymptotic normality of the adaptive mean function estimator. Our estimators are new in the context of functional time series. In Section \ref{sec:emp_study} we present a few results from an extensive simulation study and a real data analysis. The proofs of the main results are presented in the Appendix. The proofs of the  lemmas and additional technical statements are given in a Supplement \citep{hassan2024supp}. We also provide further empirical results and details on our simulation setups and the real data case study.


\section{Functional Time Series}\label{sec:FTS_model}

A functional time series (FTS) is a sequence of random functions $\{X_n\} =\{X_n(u),\ u\in I,n\in\ZZ\} \subset \Hh,$ which are temporal dependent, \emph{i.e.,} stochastically dependent with respect to the index $n$. 
Here, $I$ is a bounded domain over the real line, for instance, $I=(0,1]$. Moreover,  $\Hh=\LL^2(I)$ is the Hilbert space of real-valued, square integrable functions defined over $I$.
In applications, the index $n$ can represent the day, while $u$ can be the daily clock time rescaled to $I$. We assume that, almost surely, the paths $X_n$ belong to the Banach space $\Cc=\Cc(I)$ of continuous functions, equipped with the sup-norm $\|\cdot\|_\infty$.  


\subsection{Data}

For each $1\leq n \leq N$, the trajectory (or curve) $X_n$  is observed at the domain points $\{T_{n,i}, 1\leq i \leq  M_n\}\subset I$, with additive noise. The data points associated with $X_n$ consist of  the pairs  $(Y_{n,i} , T_{n,i} ) \in\mathbb R \times I $, where 
\begin{equation}\label{eq:data-model}
	Y_{n,i} = X_n(T_{n,i}) +  \sigma(T_{n,i})\varepsilon_{n,i},   
	\qquad 1\leq n \leq N,  \; \; 1\leq i \leq M_n.
\end{equation}
The data generating process described in \eqref{eq:data-model} satisfies the following assumptions.

\begin{assumptionH}
	\item   \label{H:stationarity} \!The series  $\{X_n\}$ is a (strictly) stationary $\mathcal H-$valued series.
	
	\item\label{H:Mn}  \!The $M_1, \dotsc,M_N$ are random  draws of an integer  variable $M\geq 2$, with expectation $\lambda$.  
	
	\item\label{H:Tni}  \!Either all the $T_{n,i}$  are independent copies of a  variable $T\in I$ which admits a strictly positive density $g$ over $I$ (independent design case), or the $T_{n,i}$, $1\leq i \leq \lambda=M_n$, are the points of the same equidistant grid of $\lambda$ points in $I$ (common design case).
	
	\item\label{H:epsni} \!The $\varepsilon_{n,i} $  are independent copies of a centered error variable $\e$ with   unit variance, and  $\sigma^2(\cdot)$ is a Lipschitz continuous function. 
	
	\item\label{H:ind} \!The series $\{X_n\}$ and the copies of  $M$, $T$ and $\varepsilon$ are mutually independent.
\end{assumptionH}

In the following, $X$ denotes  a generic random function having the stationary distribution of $\{X_n\}$. The distribution of the variable $M$ depends on $N$, namely its expectation $\lambda$   is allowed to increase with $N$. Thus, for our non-asymptotic results, the domain points $T_{n,i}$, $1\leq i\leq M_n$, $1\leq n\leq N$ are a triangular array of points. They are either obtained as random copies of $T\in I$, or they are the elements of a grid of length $\lambda$, which we consider to be equidistant for simplicity. Assumption \assrefH{H:epsni} allows for heteroscedastic errors. 

We  study the \textit{local regularity} of $X$, and thus that of the stationary distribution of $\{X_n\}$. Before providing the formal definition of the local regularity, we provide insight into this notion proposed by \citet{golovkine2022learning} in the case where the sample paths $X_n$ are not almost surely differentiable.  Let us assume that a constant $\beta>0$ exists and for any $t\in I$, $H_t\in  (0,1]$ and  $L_t \in (0,\infty)$ exist such that   
\begin{equation}\label{eq:def_lr1x}
	\EE\left[ \left\{X(u) - X(v)\right\}^2 \right] = L^2_t |u-v|^{2H_t} \{1+O(|u-v|^{\beta})\},
\end{equation}
when $u\leq t\leq v$ lie in a small neighborhood of $t$. $H_t$ is then the local Hölder exponent while $L_t$ is the local Hölder constant. They are both allowed to depend on $t$ in order to allow for curves with general patterns. Examples of processes satisfying \eqref{eq:def_lr1x} include, but are not limited to stationary or stationary increment processes \citep[][]{golovkine2022learning}. The class of multifractional Brownian motion processes with domain deformation is another example \citep[][]{wei2023adaptive}. By Kolmogorov's criterion \citep[][Theorem 2.1]{yor}, the local regularity of the process $X$ is linked to the regularity of the sample paths. Finally, the notion of local regularity extends to the case where the sample paths of $X$ admit derivatives. Condition \eqref{eq:def_lr1x} is then considered with the highest integer order derivative at $u$ and $v$ in place of $X(u)$ and $X(v)$, respectively.

\subsection{The local regularity}\label{sec_def_locr}

For any $d\in\NN$,
$\nabla ^d$ denotes  the $d-$order derivative operator, and $\mathbb R^d _+$ is the set of vectors in $\mathbb R^d$ with positive components. Let $J\subset I $ be an open interval.

\begin{assumptionH}
	
	\item\label{D:smooth2}
	For some $\delta\in\NN$, the stationary distribution of $\{X_n\}$ satisfies the following conditions~: Lipschitz continuous functions $H_\delta:J \rightarrow (0,1]$ and  $L_d:J \rightarrow (0,\infty)$, $d\in\{0,\ldots,\delta\}$, and constants $\beta_\delta >0$ and $S_\delta >0$ exist such that:  
	\begin{enumerate}
		\item with probability 1, for any $d\in\{0,\dotsc,\delta\}$, the function $\nabla^d X$  exists over $ J$, and 
		\begin{equation}\label{eq:D:L2-norm}
			0 < \underline{a}_d := \inf_{u\in J} \EE\left[\left(\nabla^d X(u)\right)^2\right] 
			\leq \sup_{u\in J} \EE\left[\left(\nabla^d X(u)\right)^2\right] =: \overline{a}_d < \infty;   
		\end{equation}
		\item there exists $\Delta_{\delta,0}>0$  such that  $\forall t,u,v\in J$ with $t-\Delta_{\delta,0}/2 \leq u \leq t \leq v \leq t+\Delta_{\delta,0}/2$, 
		\begin{equation}\label{eq:D:holder}
			\left|
			\EE\left[\left\{\nabla^\delta X(u) - \nabla^\delta X(v)\right\}^{2}\right] 
			-  L_{\delta,t}^2|u-v|^{2H_{\delta, t}}
			\right|
			\leq
			S_\delta ^2|u-v|^{2H_{\delta,t}+2\beta_\delta}.
		\end{equation}
	\end{enumerate}
	
\end{assumptionH}

\begin{definition}\label{def:space_family}
	\!	Let $\Xx(\delta+H_\delta , \boldsymbol L_\delta  ;J)$ denote the class of stochastic processes $X$ with continuous paths satisfying \assrefH{D:smooth2}, where $\boldsymbol L _\delta = (L_0,\ldots,L_\delta)\in \mathbb R^{\delta +1} _+$, and 
	\begin{equation*}
		0<\inf_{u\in J} H_{\delta, u}  \leq  \max_{u\in J} H_{\delta,u} <1 \quad \text{and} \quad 0<\min_{0\leq d\leq \delta}\inf_{u\in J} L_{d,u}  \leq \max_{0\leq d\leq \delta} \sup_{u\in J} L_{d,u} < \infty.
	\end{equation*}
\end{definition}

The following result describes the embedding structure of the spaces $\Xx(\delta+H_\delta , \boldsymbol L_\delta  ;J)$. The proof  is given in \citet{hassan2024supp}. 

\begin{lemma}\label{lem:regularity-1} 
	Assume that $X$ belongs to $\mathcal{X}(\delta+H_\delta,\boldsymbol L_\delta,J)$  for some $\delta\in\NN^*$,  $J$ an open sub-interval of $I$,   $0 < H_\delta < 1$,  and a bounded vector-valued function $\boldsymbol L_\delta\in \mathbb R^{\delta+1}  _+$. 
	Then, for any $d\in\{0,\dotsc,\delta-1\}$, $X$ belongs to $\mathcal{X}(d+H_d,\boldsymbol L_d, J)$ with $H_d \equiv 1$ and some bounded vector-valued function $\boldsymbol L_d\in \mathbb R^{d+1} _+$.
\end{lemma}

The parameters defining the local regularity are formally defined  in the following. 

\begin{definition}\label{def:par_family}
	If $X \in \Xx(\delta+H_\delta , \boldsymbol L_\delta  ;J)$, with $\delta \in \NN$ and  $0 <H_{\delta,t} <1 $,   the \emph{local regularity} of $X$ at $t$, an interior point of $I$, is defined by the parameters 
	$	\alpha_t = \delta + H_{\delta,t}$ and $ L_{\delta,t}^2$.
\end{definition}

For the purposes of the applications we have in mind, when $\delta \geq 1$, estimating the Hölder constants $ L_{d,t}$ for $0\leq d \leq \delta-1$ is worthwhile, and thus will not be considered. 
We next present examples of FTS and their regularity parameters. 

\begin{example}
	\label{ex:WN}\normalfont 
	Let $\{\xi_n\} $ be a sequence of i.i.d. multifractional Brownian motion (MfBm) of Hurst exponent function $H_\xi:\mathbb R_+  \rightarrow (0,1)$. That means $\xi_n$ are independent copies of $\xi$, a centered Gaussian process with covariance function 
	\begin{equation*}
		\EE\left[\xi(u)\xi(v)\right] =  D(H_{\xi,u},H_{\xi,v} )\left[ u^{H_{\xi,u}+H_{\xi,v}} +  v^{H_{\xi,u}+H_{\xi,v}} - |v-u|^{H_{\xi,u}+H_{\xi,v}}\right] , \qquad u, v\geq 0,
	\end{equation*}
	where
	\begin{equation*}
		D(x,y )=\frac{\sqrt{\Gamma (2x+1)\Gamma (2y+1)\sin(\pi x)\sin(\pi y)}} {2\Gamma (x+y+1)\sin(\pi(x+y)/2)} , \qquad D(x,x) = 1/2, \qquad x,y >0.
	\end{equation*}	
	See, \emph{e.g.}, \citet{balanca2015}  for the formal definition of the MfBm.  The  fractional Brownian motion is an  MfBm with constant Hurst index function. 
	For any bounded interval $I \subset \mathbb R_+$, it can be shown that $\xi \in \Xx(H_\xi , 1  ;I)$ provided $H_\xi$ is twice continuously differentiable \citep{wei2023adaptive}. 
	Note that defining
	$
	\eta (t) = \int_a^t \xi(u) du$, $t \geq 0$,
	for some $a \geq 0$, we have $\eta \in \Xx(1+H_\xi , \boldsymbol L_1  ;I)$, where $\boldsymbol L_1(t) = (\operatorname{Var}(\xi(t)),1 )$. Repeatedly applying the integral operator yields examples of processes with any non-integer  $\alpha_t>1$ in Definition \ref{def:par_family}. 
\end{example}

\begin{example}[FAR(1) model]\label{ex:far-1}\normalfont 
	Let $\{X_n\} $ be the zero-mean, stationary Functional AutoRegressive (FAR) time series that is the stationary solution of the equation
	\begin{equation}\label{eq:far-1}
		X_n(t) = \Psi(X_{n-1})(t) + \xi_n(t),
		\qquad t \in I\subset \mathbb R_+,\quad n\in\ZZ,
	\end{equation}
	where $\{\xi_n\} $ is an MfBm functional white noise as in Example \ref{ex:WN}, with the twice continuously differentiable Hurst index function $H_\xi \in (0,1)$. We next assume that  $\Psi$ is the integral operator
	\begin{equation*}
		\forall x \in \Cc,\quad 	\Psi(x)(t) = \int_I \psi(s,t) x(s) ds, \qquad  \text{ with } \quad \iint_{I\times I} \psi^2(s,t)dsdt < 1.
	\end{equation*}
	The stationary solution for \eqref{eq:far-1} then exists \citep[see, for instance,][Section 8.8]{kokoszka2017introduction}.
	Assume further that constants $C>0$, $H_\psi \in (0,1]$ exist such that 
	\begin{equation*}
		\sup_{u\in I}H_{\xi,u} < H_\psi  \leq 1 \qquad \text{and} \qquad 	\left| \psi(s,u) - \psi(s,v) \right|^2 \leq C|u-v|^{2 H_\psi}, 
		\quad \forall s,u,v\in I .
	\end{equation*}
	Then,  $\{X_n\} $ belongs to  $\Xx(H_\xi , 1  ;I)$. See \citet{hassan2024supp} for the details. 
\end{example}

\subsection{Weak dependence}

We consider a general notion of weak dependence which allows for a refined study of the local regularity of FTS. More precisely, we reconsider the concept of \textit{$\LL^p_\Hh-m-$approximability}  \cite[see][]{hormann2010weakly} in the context of $(\Cc,\|\cdot\|_\infty)$-valued (instead of $\Hh$-valued) random processes. In this way, the type of weak dependence between the curves $X_n$ is inherited by the sequences $\{X_n(t)\}$, for all $t\in I$. The general idea with the weak dependence type considered by Hörmann and Kokoszka is to approximate $\{X_n\}$ by an $m$-dependent sequence $\{X_n^{(m)},m\geq 1\}$ such that, for every $n\in\ZZ$, the sequence $\{X_n^{(m)}, m\geq 1 \}$ converges in some sense to $X_n$ as $m\to\infty$.
The limiting behavior of the original process can then be obtained from that of its coupled $m$-dependent sequences provided they are sufficiently close to the original process. 

Some more notations are needed~: $\left<\cdot,\cdot\right>_\Hh$ and $\|\cdot\|_\Hh$ denote the inner product of the Hilbert space $\Hh$ and the associated norm respectively. 
For $p\geq 1,$   $\LL^p$ is the space of real-valued variables $Z$ with $\nu_p(Z) = \left(\EE\left[ |Z|^p \right]\right)^{1/p} < \infty$.
Moreover,  $\LL^p_\Hh$ and $\LL^p_\Cc$ are the spaces of $\Hh$-valued and  $\Cc$-valued random functions $X$ with $\nu_p\left(\|X\|_\Hh \right) < \infty$ and  $\nu_p\left(\|X\|_\infty \right) < \infty$ respectively.  

\begin{definition}\label{H:Lpm-approx} 
	The stationary FTS
	$\{X_n\}$ is $\LL^p_\Cc-m-$approximable with $p\geq 1$ if~:
	\begin{enumerate}
		\item  $\{X_n\}\subset \LL^p_\Cc$ admits a moving average (MA) representation, \emph{i.e.,}
		\begin{equation}\label{eq:MArepresentation}
			X_n = f(\xi_n,\xi_{n-1},\ldots)
		\end{equation}
		with  $\{\xi_n\}$  independent copies of $\xi\in S$, $S$ a measurable space and  $f\!:\!S^\infty\!\to\Cc$ 
		measurable. 
		\item\label{cndt_3} For every $n\in\ZZ$, let $\{\xi_k^{(n)}, k\in\ZZ \}$ be a sequence of independent   copies of $\xi$  defined over the same probability space.
		The coupled version of $X_n$ is defined by
		\begin{equation*}
			X_n^{(m)} = f(\xi_n,\xi_{n-1},\ldots,\xi_{n-m+1},\xi_{n-m}^{(n)},\xi_{n-m-1}^{(n)},\ldots).  
		\end{equation*}
		\item\label{cndt_4} The sequence $\{X_n^{(m)}, m\geq 1\}$ converges to $X_n$ as $m\to\infty$ in the sense that 
		$$
		\sum_{m\geq 0} \nu_p\left( \|X_m - X_m^{(m)}\|_\infty \right) <\infty.
		$$ 
	\end{enumerate}    
\end{definition}

As in \cite{hormann2010weakly}, having $p\geq 4$ will be convenient for our applications.
Moreover, in Assumption \assrefH{H:Lpmapprox} below we impose stronger restrictions on the rate of convergence of the coupled sequences, see also  \cite{rice2019inference}.

\begin{assumptionH}
	\item\label{H:Lpmapprox} The stationary FTS  $\{X_n\}$ is $\LL^p_\Cc-m-$approximable with  some  $p\geq 4$ such that constants $C>0$ and $\alpha >3/2$ exist and $\nu_p ( \|X_m - X_m^{(m)}\|_\infty  )\leq C m^{-\alpha}$, $m\geq 1$.
\end{assumptionH}

The basic  properties of $\LL^p_\Hh-m-$approximability established by \citet[Lemma 2.1]{hormann2010weakly} remain true with Definition \ref{H:Lpm-approx}, (see our Lemma \ref{lem:lp-m-approx-1} in the Appendix). Moreover, Lemma \ref{lem:lp-m-approx-2}  shows that $\LL^p_\Cc-m-$approximability entails the pointwise $\LL^p-m-$approximability of $\{X_n(t)\}$, for all $t\in I$.
Note also that $\LL^p_\Cc-m-$approximability implies $\LL^p_\Hh-m-$approximability, because the $\|\cdot\|_\Hh$  is bounded by the sup-norm. More generally, Definition \ref{H:Lpm-approx} can be considered with other Banach spaces $\mathcal C$  than $(\mathcal C(I),\|\cdot\|_\infty)$.
For instance, when $\mathcal C$ is the real line, our definition of $\LL^p_\Cc-m-$approximability becomes the $\LL^p-m-$approximability for scalar times series, see \citet{wu2005nonlinear}.

We now show that some common FTS models are $\LL^p_\Cc-m-$approximable, in the sense of Definition \ref{H:Lpm-approx}. 
Let $\Ll = \Ll(\Cc, \Cc)$ denote the space of bounded linear operators on $(\Cc,\|\cdot\|_\infty)$.
For any Hilbert-Schmidt operator $A$, let 
$\vertiii{A}_\infty = \sup \{\| Ax \|_\infty : \| x \|_\infty \leq 1\}$.
The 
justification for the following examples is given in \citet{hassan2024supp}. 

\begin{example}[FAR(1) model revisited]\label{ex:far} \normalfont
	Consider the model in Example \ref{ex:far-1}, with $\Psi\in \Ll$   such that $\vertiii{\Psi}_\infty < 1$, and a zero-mean i.i.d. sequence $\{\xi_n\} \subset \LL_{\Cc}^p$. 
	By \citet[Theorem 3.1]{bosq2000linear}, there exists then, a unique mean zero,  stationary solution $\{X_n\}\subset\Cc$ of the FAR(1) equation \eqref{eq:far-1}, provided $p\geq 2$. 
	Then, $\{X_n\}$ is $\LL^p_\Cc-m-$approximable.
\end{example}

\begin{example}[Functional linear process]\label{ex:linear} \normalfont
	Let $\{X_n\}$ be the \textit{linear process}  
	defined as  
		$X_n = \sum_{j=0}^{\infty}\Psi_j(\xi_{n-j})$,  with $\{\xi_j\}\subset\LL^p_\Cc$  i.i.d., $\EE (\xi_j) = 0$,
	and   the  operators
	$\Psi_j \in \Ll$ 	satisfy 
		$ 
		\sum_{j=1}^{\infty} j\vertiii{\Psi_j}_\infty< \infty$. 
	Then, $\{X_n\}$ is $\LL^p_\Cc-m-$approximable. 
\end{example}

\begin{example}[Product Model]\label{ex:product}\normalfont
	Suppose that $\{Y_n\} \subset \LL_\Cc^p$ and $\{U_n\} \subset \LL^p$ are two independent $\LL^p-m-$approximable sequences. Their  MA representations are $Y_n = g_Y(\eta_1, \eta_2, \dots\!)$ and $U_n = g_U(\gamma_1, \gamma_2, \dots)$, where $\{\eta_n\}$ and $\{\gamma_n\}$ are two i.i.d. random sequences. 
	The sequence $\{X_n\} \subset \LL_\Cc^p$ with $X_n(t) = U_n Y_n(t)$, $t\in I$,  is then $\LL_\Cc^p-m-$approximable sequence with the  i.i.d. variables $\xi_n = (\eta_n, \gamma_n)$ in the MA representation \eqref{eq:MArepresentation}. 
\end{example}

\begin{example}[Functional ARCH(1)]\label{ex-FARCH}\normalfont
	Let $c(\cdot) \in \Cc$ be a positive function and $\{\xi_n\}$ 
	independent copies of $\xi\in \LL_\Cc^p$. Let $\beta(s,t)$ be a continuous, non-negative 
	function. 
	Then,
	\begin{equation}\label{eq:arch-1}
		Y_n(t) = \xi_n(t)\sigma_n(t) \quad \text{ with } \quad 
		\sigma_n^2(t) = c(t) + \int_I \beta(s,t) Y_{n-1}^2(s) ds, \quad t\in I,n\in\ZZ,
	\end{equation}
	is the  functional AutoRegressive Conditional Heteroskedastic (ARCH) series of order 1. If  
	\begin{equation*}
		\text{for some $p > 0$,}\quad 	\EE[  H^{p/2} (\xi^2)]  < 1\quad\text{with }\quad H(\xi^2) = \underset{t \in I}{\sup} \int_0^1 \beta(s,t)\xi^2(s) ds,
	\end{equation*}
	then  \eqref{eq:arch-1} 
	has a unique, strictly stationary solution $\{Y_n\}$ \citep[see][Theorem 2.2]{hormann2013arch}. Moreover, the solution  is $\LL^p-m-$approximable. 
\end{example}


\section{Estimation of the local regularity parameters}\label{sec:local-reg-estimator}

For simplicity, we first consider the case  where the sample paths of $X$ are almost surely non-differentiable,  which means $\delta=0$ in the definitions  in Section \ref{sec_def_locr}. Functional data from fields such as energy, environment,  chemistry and physics, medicine, meteorology, are  often very irregular and thus considering non-differentiable curves $X_n$ seems realistic. The case $\delta \geq 1$ is discussed 
in Section \ref{sec_reg_diff}.

\subsection{The case of non-differentiable sample paths} \label{sec:local-reg-estimator_nd}
Set $t\in J$. 
We simplify the notation and denote by  $(H_t,L_t^2)$, instead of $(H_{0,t},L_{0,t}^2)$, the local regularity parameter at point $t$.
Let $\Delta \leq \Delta_{0,0}$ and  $t_1,t_2, t_3\in J$  such that $t_3 - t_1 = \Delta$ and $t_2=t = (t_1+t_3)/2$.  
Using the definition of the local regularity, the following proxy values of $H_t$  and $L_t^2$ are considered,
\begin{align}
	\widetilde H_t  = \widetilde{ H }_t (\Delta) &= \frac{\log(\theta(t_1,t_3)) - \log(\theta(t_1,t_2))}{2\log(2)}, \label{eq:proxyH}
	\\ 
	\widetilde L_t^2  =  \widetilde{L}_{t}^2(\Delta) &=  \theta\left(t_1, t_3\right) \Delta^{-2\widetilde H_t   }, \qquad \text{where } \qquad \theta(u,v)=\EE[ \left\{X(u) - X(v)\right\}^2 ] .
\end{align}
Lemma \ref{lem:convergence_proxy} in Appendix states that $\widetilde H_t$ and $\widetilde{L}_{t}^2 $  converge  to $H_t$ and $L_t^2$  as $\Delta\rightarrow 0$.
Moreover, we have $\widetilde{L}_{t}^2 =L_t^2$ and $\widetilde H_t =H_t$ if $S_0^2=0$ in \eqref{eq:D:holder}.  
Our estimators of $H_t$ and $L_t^2$ are obtained by plugging  the estimators of $\theta(\cdot,\cdot)$ into the definition of  $\widetilde H_t$ and $\widetilde{L}_{t}^2$, respectively.

\paragraph*{Presmoothing step}
The estimation of $\theta(u,v)$ implies the reconstruction of the  curves $X_1,\ldots,X_N$ at the points $u$ and $v$, using data as described in  \eqref{eq:data-model}.
To preserve the stationarity of the reconstructed curves,  we use the same linear presmoothing estimator for all $X_n$. Given the sample points $(\Yni , \Tni )$, $1 \leq i \leq M_n$, the presmoothing estimator of $X_n$ is defined as
\begin{equation}\label{eq:pre-smooth}
	\widetilde{X}_n(u) = \sum_{i=1}^{M_n} W_{n,i}(u) Y_{n,i},\qquad u\in J,\ n=1,\ldots,N,
\end{equation}
where the weights $\{W_{n,i}\}_{i=1\ldots M_n}$ depend on $(M_n,T_{n,1},\ldots,T_{n,M_n})$ and some presmoothing parameter. We impose the following assumptions on the presmoothing estimator.

\begin{assumptionH}
	
	\item\label{H:erreur_oublie} The error variable $\varepsilon$ from \assrefH{H:epsni} has finite moment of order $p$ with $p$ from \assrefH{H:Lpmapprox}.

	\item\label{H:pre-smooth:weights} The sums of the absolute values of the weights $W_{n,i}(u) $ are bounded by a constant.

	\item\label{H:pre-smooth:convergence} Constants $B,\tau>0$ exist such that 
	$
	R_2(\lambda)  = \sup_{u\in J} \EE[\{\widetilde{X}_n(u) \!- \!X_n(u)\}^2 ] \leq B \lambda^{-\tau}.
	$
\end{assumptionH}

Assumption \assrefH{H:pre-smooth:weights} is always satisfied with the constant equal to 1 when the weights $W_{n,i}$ are non-negative, and this is the case for the Nadaraya-Watson estimator. Assumption \assrefH{H:pre-smooth:convergence} is a mild condition. For instance, when $T$ admits a density and \assrefH{H:boundMn} below holds true, \assrefH{H:pre-smooth:convergence} can be guaranteed when the density of $T$ 
stays away from zero, and that density together with the sample paths of $X$, satisfy some mild smoothness assumptions \citep[e.g.,][]{tsybakov2009}. 

\paragraph*{Local regularity estimators}
Given the presmoothed curves $\widetilde{X}_n$,    the estimator of $\theta(u,v)$ is
\begin{equation}\label{def_theta_hat}
	\widehat\theta(u,v) = \frac1{N} \sum_{n=1}^{N} \left(\widetilde{X}_{n}(v) - \widetilde{X}_{n}(u)\right)^{2}, \quad u,v\in J.
\end{equation}
Our estimators of $H_t$ and $L_t^2$ are then defined as,
\begin{equation}\label{eq:H0hat}
	\widehat H_t = \frac{\log(\widehat\theta(t_1,t_3)) - \log(\widehat\theta(t_1,t_2))}{2\log(2)}\quad \text{ and  } \quad 
	\widehat{L}_{t}^2 = \frac{\widehat\theta\left(t_1, t_3\right)}{\Delta^{2\widehat H_t}}.
\end{equation}

\begin{theorem}\label{thm:regularity:H}
	Assume that \assrefH{H:stationarity} -- \assrefH{H:pre-smooth:convergence} hold true, and $\widetilde{ H }_t, \widetilde{ L }^2_t$ are defined with $\Delta \leq \Delta_{0,0}$. A constant $C>0$
	exists such that, for any  $\varphi\in(0,1)$ satisfying the conditions
	\begin{align}
		\Delta^{2\beta_0} S_0^2 &< \frac{L_t^2\log(2)}{4}\varphi,\label{eq:thmH:proxy}\\
		\lambda^{-\tau/2} &< C L_t^2 \varphi\Delta^{2H_t},\label{eq:thmH:tau-phi}
	\end{align}
	we have 
	\begin{equation*}
		\PP\left( \left|\widehat H_t - H_t \right| > \varphi \right)
		\leq
		\frac{ \mathfrak{f}_0 }{ N \varphi^2 \Delta^{4 H_t} }\\
		+ \mathfrak{b} \exp{ \left( - \mathfrak{g}_0 N \varphi^2 \Delta^{4 H_t} \right) },
	\end{equation*}
	for some universal constant $\mathfrak{b}$, provided $\lambda$ is sufficiently large.  The constant  $C$ depends on  $\overline{a}_0$ from \eqref{eq:D:L2-norm} and  $B$ from \assrefH{H:pre-smooth:convergence}, and $\mathfrak{f}_0$ and $\mathfrak{g}_0$ are determined by the dependence structure of $X$.
\end{theorem}

\begin{theorem}\label{thm:regularity:L}
	Assume that the conditions of Theorem \ref{thm:regularity:H} hold true.
	A constant $\widetilde{C}>0$ exists, 
	such that for any  $\varphi,\psi\in(0,1)$ satisfying the additional conditions
	\begin{align}
		3\Delta^{-2\varphi} \Delta^{2\beta_0} S_0^2 &< \psi,\label{eq:thmL:proxy}\\
		6 L_t^2 \Delta^{-2\varphi} \varphi |\log\Delta| &< \psi ,\label{eq:thmL:phi-psi}\\
		\lambda^{-\tau/2} &< \widetilde{C}\Delta^{2\varphi} \psi\Delta^{2H_t} ,\label{eq:thmL:tau-psi}
	\end{align}
	we have 
	\begin{multline*}
		\PP\left( \left| \widehat{L}_{t}^2 - L_{t}^2 \right| > \psi \right) 
		\leq 
		\frac{\mathfrak{c}_0}{N\psi^2 \Delta^{4H_t + 4\varphi}} 
		+ \frac{ \mathfrak{f}_0 }{ N \varphi^2 \Delta^{4 H_t} } \\
		+\mathfrak{b}\exp{ \left( - \mathfrak{l}_0 N\psi^2 \Delta^{4H_t + 4\varphi} \right) }
		+ 4\mathfrak{b} \exp{ \left( - \mathfrak{g}_0 N \varphi^2 \Delta^{4 H_t} \right) } ,
	\end{multline*}
	for some universal constant $\mathfrak{b}$, provided $\lambda$ is sufficiently large. The constant  $\widetilde{C}$ depends on $\overline{a}_0$ and $B$, while the constants  $\mathfrak{c}_0$, $\mathfrak{f}_0$, $\mathfrak{g}_0$, $\mathfrak{l}_0$ are determined by the  dependence structure of $X$.
\end{theorem}

\begin{remark}
	Since $\Delta<1$ and $\varphi>0,$ the condition 
	\eqref{eq:thmL:phi-psi} 
	implies $6 L_t^2 \varphi |\log\Delta| < \psi$. Thus, when  $\psi$ decreases to $0$, $\varphi |\log\Delta|$ and $\Delta^{-2\varphi}$ converge to $0$ and $1$, respectively.
\end{remark}

\begin{remark}
	The role of \eqref{eq:thmH:proxy} and \eqref{eq:thmL:proxy} is to control the bias between the  parameters and their proxies $(\widetilde{H}_t,\widetilde{L}_t^2)$. When $S_0=0$, the convergence rates of $\widehat{H}_t$ and $\widehat{L}_t^2$ are given by $\varphi =\Oo(\lambda^{-\tau/2}\Delta^{-2H_t})$ and $\psi =\Oo(\varphi(\lambda)|\log(\Delta)| ) = \Oo(\lambda^{-\tau/2}\Delta^{-2H_t} |\log(\Delta)| ), $  respectively.
	Meanwhile, when $S_0\neq 0$, the conditions \eqref{eq:thmH:proxy} and \eqref{eq:thmL:proxy} tend to decrease these rates of convergence, that are $\varphi =\Oo(\lambda^{-\tau\beta_0/(2\beta_0+2H_t)})$ and $\psi =\Oo(\lambda^{- \tau\beta_0/(2\beta_0+2H_t)} |\log(\Delta)|)$  under the condition that $\Delta=\Oo(\lambda^{-\tau/(4\beta_0+4H_t)})$.
\end{remark}

In applications where the local regularity estimation serves some specific purposes such as adaptive estimation of the mean, and (auto-)covariance functions, it will suffice to consider $\varphi = \left(\log\lambda\right)^{-2}$ and 
$\psi  =   \left(\log\lambda\right)^{-1}$. The only choice the statistician has to make is that of $\Delta$, for which we propose 
$\Delta = \exp\left( -(\log\lambda)^\gamma\right)$ for some $\gamma\in(0,1)$. See Section~\ref{sec:emp_study:loc_reg}.

\subsection{Regularity estimation for differentiable paths}\label{sec_reg_diff}

Following \cite{golovkine2022learning}, we now construct an estimator of the local regularity when  $X$ restricted to $J$, belongs to $\mathcal{X}(\delta + H_\delta, \boldsymbol{L}_\delta; J)$ with $\delta \in \NN$.
Lemma~\ref{lem:regularity-1} indicates that when $\delta\geq 1,$ $\nabla^d X$ restricted to $J$ belongs to the class $\mathcal{X}(H_d, \boldsymbol{L}_d;J)$ with $H_d=1$ if $d<\delta$.
With at hand an estimator $\widehat \lambda$ of $\lambda$ and a suitable  decreasing function  $ \varphi(\lambda)$, for instance $\varphi(\lambda) =  (\log\lambda)^{-2}$, the  Theorem~\ref{thm:regularity:H} suggests as estimator of $\delta$ the nonnegative integer
\begin{equation*}
	\widehat\delta = \min  \{ d\in\NN : \widehat H_{d,t}  < 1 - \varphi(\widehat \lambda)\}, 
\end{equation*}
where  $\widehat{H}_{d,t}$ is an estimator of the local regularity exponent parameter of $\{\nabla^d X_n\}$ at $t$.
More precisely, for $d\geq 1$, given a presmoothing  estimator $\widetilde{\nabla^d X}_n(u)$ of ${\nabla^d X}_n(u)$, for $u\in J$,  the estimators of $H_{d,t}$ and $L^2_{d,t}$ are
\begin{align*}
	\widehat{H}_{d,t} &= \frac{\log \widehat\theta_d(t_1,t_3) - \log \widehat\theta_d(t_1,t_2)}{2\log (2)}, \\ 
	\widehat{L}_{d,t}^2 &= \frac{\widehat{\theta}_d\left(t_1, t_3\right)}{\Delta^{2\widehat{H}_{d,t}}} \qquad \text{ where } \qquad \widehat\theta_d(u,v) = \frac1{N} \sum_{n=1}^{N}
	\left(
	\widetilde{\nabla^d X}_n(u) - \widetilde{\nabla^d X}_n(v)
	\right)^2. 
\end{align*}
A natural estimator of the local regularity parameter $\alpha_t$ is then 
	$\widehat\alpha_t = \widehat\delta + \widehat H_{\widehat\delta,t}$.
The procedure is summarized in 	Algorithm \ref{algo:alpha_dc}.
A detailed justification and the concentration bounds for these estimators can be found in \citet{hassan2024supp}.

\vspace{-0.3cm}	
\begin{algorithm}
	\DontPrintSemicolon 
	\KwIn{Function $\varphi(\lambda)$; integers $M_1,\ldots,M_N$; data points $(Y_{n,i} , T_{n,i} )$ generated  as in  \eqref{eq:data-model}, $1\leq i \leq M_n$, $1\leq n\leq N$}
	\KwOut{Estimation of $\alpha_t = \delta + H_{\delta,t}$}
	$\widehat{\lambda} \gets N^{-1}(M_1+\ldots+M_N)$,
	$d \gets 0$\;
	Compute $\widehat{H}_{0,t}$ as in  \eqref{eq:H0hat}\;
	\While{$\widehat{H}_{d,t}  \geq  1 - \varphi(\widehat{\lambda})$}{
		Estimate the $(d+1)$-th derivative of the trajectories of $\{X_n\}$\;
		Calculate $\widehat{H}_d$ using the estimated trajectories of the $(d+1)$-th derivatives\;
		Set $d \gets d+1$\;
	}
	\Return{$d+\widehat{H}_{d,t}$  }\;
	\caption{Estimation of the local regularity $\alpha_t$ with differentiable sample paths}
\label{algo:alpha_dc}
\end{algorithm}
\vspace{-0.2cm}

\vspace{-0.4cm}
\section{Adaptive mean and autocovariance functions estimators}\label{sec:adapt_mean_acov}

We consider a stationary FTS $\{X_n\}\subset \LL^p_\Cc$, with $p\geq 4$, defined over $I=(0,1]$. The mean function  is 
$
\mu(t)=\EE[X_n(t)] $, $t\in I,
$
and for $\ell\in\NN^*$, its lag-$\ell$ cross-product and lag-$\ell$ autocovariance functions are 
$$
\gamma_\ell(s,t) = \EE\left[X_{\ell}(s) X_{n+\ell}(t)\right] \quad \text{and} \quad 	\Gamma_\ell(s,t) = \gamma_\ell(s,t) - \mu(s)\mu(t) ,\qquad s,t\in I,
$$
respectively.  In this section we propose nonparametric estimates of $\mu(t)$ and $\gamma_\ell(s,t)$. Our estimates adapt to the regularity of $X$ and to the best of our knowledge, are the first of this kind in the context of weakly dependent FTS. For simplicity, we assume that the sample paths $X_n$ are not differentiable ($\delta=0$), and we simply denote by $(H_t,L_t^2)$ the local regularity parameters at point $t$. 
Let  $(\widehat H_t,\widehat L_t^2)$ be the estimators of $(H_t,L_t^2)$ defined according to \eqref{eq:H0hat}.

In the independent design case, let $\widehat\lambda = N^{-1} \sum_{n=1}^N M_n$ be the empirical mean of the number of observation times per curve. In both the independent and common design cases, let $\widehat \sigma^2(t)$ be a  consistent estimator of errors' variance $\sigma^2(t)$ at $t$. 
A simple choice is
\begin{equation}\label{eq:sigma_estimator}
\widehat \sigma^2 (t):=  
\frac{1}{N}\sum_{n=1}^N \frac{1}{2} \left(Y_{n,i(t)} - Y_{n,i(t)+1}\right)^2,
\end{equation}
where, for each $n$, $i(t), i(t)+1$ are the indices of the two closest domain points $\Tni$ to $t$. For each $1\leq n \leq N$, we consider the Nadaraya-Watson (NW) estimator of the trajectory $X_n,$
\begin{equation}\label{LP_est_v}
\widehat X_n(t,h) = \sum_{i=1}^{M_n}W_{n,i}(t;h)  \Yni , \qquad 
\; W_{n,i}(t;h)  =  K\!\left( \!\frac{ \Tni -t}{h}\! \right)\!\!\left[\sum_{k=1}^{M_n} K\!\left( \!\frac{T_{n,k} - t}{h}\! \right)\!\right]^{-1}\!\!\!\!, 
\end{equation}
where $h$ is the bandwidth parameter considered in some range $\Hh_N$,  $K$ is a non-negative, symmetric  and bounded kernel with the support in $[-1,1]$, and the convention $0/0=0$ applies. With at hand the estimates $\widehat X_n(t,h)$, we follow the  `\textit{smooth first, then estimate}' approach and define the estimators for the mean  and the lag-$\ell$ cross-product functions under the form of empirical estimators with the true values of the curves replaced by the smoothed ones.

Before formally defining our adaptive estimators, let us point out that the estimator $\widehat X_n(t;h)$ is degenerate if there are no domain points $ \Tni$ in $[t-h, t+h]$. This situation occurs with any smoothing-based method, and is more likely when the curves are sparsely sampled in their domain. When a curve is not observed in the neighborhood of $t$, this means that it does not carry useful information about $X_n(t)$, and should thus be dropped from the data set when estimating $\mu(t)$ or $\gamma_\ell(s,t)$.
The neighborhood of $t$ is defined by  $h$. A trade-off must be found between the large bias induced by large $h$ and  the  large variance resulting from dropping curves when $h$ is small. For this purpose, let $\mathds{1}\{\cdot\}$ denote the indicator function and    
\begin{equation}\label{eq:def_pi}
\pi_n (t;h) = 1 \quad \text{if } \; \sum_{i=1}^{M_n} \mathds{1}\{\lvert  \Tni-t\rvert\leq h\} \geq 1, \quad \text{ and } \; \pi_n (t;h) = 0  \text{ otherwise}.
\end{equation}
The number of curves $X_n$ with at least one observation in the interval $[t-h,t+h]$ is then
\begin{equation*}
P_N(t;h) =\sum_{n=1}^N \pi_n (t;h).
\end{equation*}
We finally denote the conditional expectation given the $M_n$ and the realizations of $T$ by~: 
$$
\EEMT(\cdot) = \EE \left(\cdot \mid \! M_n, \mathcal \{\Tni, 1\leq i \leq M_n\}, 1\leq n \leq N \right). 
$$

\subsection{Adaptive mean function estimator}

There are several contributions to the problem of estimating the mean function in the context of stationary FTS.
First, if the curves are fully observed without error, under the $L^2_\Hh-m-$approximation assumption, the empirical mean is a $\sqrt{N}$-consistent estimator as per \citet{hormann2012functional}.
A central limit theorem for the empirical mean is also established under cumulant mixing dependence.
\citet{sabzikar2023tempered} have defined a broad class of models for FTS that can be used to quantify near long-range dependence. They established rates of consistency for the  empirical mean function assuming error-free, fully observed sample paths. Recently, authors have been focusing on the case where FTS are discretely sampled with an additive noise.
\citet{chen2015simultaneous} and \citet{li2023statistical} propose a method, based on $B$-splines, for constructing simultaneous confidence bands for the mean function under physical-dependence and infinite average FTS models, respectively. Their procedures assume an equidistant common design and the mean function is at least continuously differentiable.
\citet{rubin2020sparsely} propose a local linear estimator of the mean function when the design is random and when the curves are sparsely observed.
They derive asymptotic results assuming a mean function that is twice differentiable and two types of dependence conditions, namely  cumulant mixing and strong mixing conditions. We here propose an adaptive nonparametric mean function estimator for irregular mean functions, and derive its asymptotic normality.

Let $t\in I$ be fixed. Our adaptive mean function pointwise estimator is 
\begin{equation}\label{eq:hatmu}
\widehat \mu_N^* (t) = \widehat \mu_N (t;h_{\mu}^*) \quad \text{ with } \quad \widehat \mu_N (t;h) = \frac{1}{P_N(t;h) }\sum_{n=1}^N  \pi_n (t;h) \widehat X_n(t;h),
\end{equation}
where $h_{\mu}^*$ is an adaptive, optimal bandwidth. To define the selection rule for $h$, let
\begin{equation} \label{eq:mu:risk-boundV}
R_\mu(t;h, H_t, L_t^2, \sigma^2(t)) = L_t^2 h^{2  H_t}\mathbb{B}(t;h,2  H_t) +  \sigma^2\! (t)  \mathbb{V}_\mu(t;h) + {\mathbb{D}}_\mu(t;h)/P_N(t;h),
\end{equation}
and $R_\mu(t;h):=	R_\mu(t;h, H_t, L_t^2, \sigma^2(t))$, where 
\begin{align*}
\mathbb{V}_\mu(t;h) &= \dfrac{1}{P^2_N(t;h)} \sum_{n=1}^N  \pi_n(t;h) c_n(t;h) \max_{1\leq i \leq M_n} \left\lvert W_{n,i}(t;h)\right\rvert,
\\
\mathbb{B} (t;h,\alpha) &= 
\frac{1}{P_N(t;h)} \sum_{n=1}^{N}\pi_n (t;h) c_n(t;h) b_n(t;h,\alpha), \text{ with } c_n(t;h) = \sum_{i=1}^{M_n} \left|W_{n,i}(t;h)\right|
,\\
\mathbb{D}_\mu (t;h) &= \EE\left[\{X_0(t) - \mu(t)\}^2\right] + 2\sum_{\ell = 1}^{N-1}p_\ell(t;h) \EE\left( \{X_0(t) - \mu(t)\} \{ X_\ell(t) - \mu(t)\} \right),\\
\text{with}&\quad p_\ell(t;h) = \sum_{i=1}^{N-\ell}\frac{\pi_i(t;h)\pi_{i+\ell}(t;h)}{P_N(t;h)},\quad
b_n(t;h,\alpha) = 
\sum_{i=1}^{M_n} \left|\frac{ \Tni - t}{h}\right|^{\alpha} \left\lvert W_{n,i}(t;h)\right\rvert.
\end{align*}
The three terms on the right-hand side of \eqref{eq:mu:risk-boundV} can be interpreted as a bias, a stochastic and a penalty term, respectively. Regarding the latter, which is specific to the FDA framework, ${\mathbb{D}}_\mu(t;h)/P_N(t;h)$ increases as $h$ decreases because more curves are excluded from the mean estimation. 
We will show that $2R_\mu(t;h)$ is a sharp bound of the quadratic risk $\EEMT\left[(\widehat \mu_N(t;h) - \mu(t))^2\right]$ over a wide  grid $\Hh_N$ of bandwidths. Note that $ c_n(t;h)\equiv 1$ in the case of NW estimator with a non-negative kernel. Note also that under the $\LL^4_\Cc-m$-approximation assumption,  the autocovariances of the time series $\{X_n(t), n \geq 1\}$ are absolutely summable \cite[see][Lemma 4.1]{hormann2010weakly}.  This means that, without using any additional information on the FTS model, we can simply take absolute values and bound $\mathbb{D}_\mu (t;h)$
by a constant equal to the limit of the series of the absolute values of the autocovariances. For now, we suppose that $\mathbb{D}_\mu (t, h)$ is given.

The bandwidth $h_\mu^*$ is selected to minimize an estimator of $R_\mu(t;h)$. 
More precisely, 
\begin{equation}\label{eq:mu:risk-minimization}
h_{\mu}^* \in \underset{h \in \Hh_N}{\argmin}\; \widehat R_\mu(t;h) \qquad \text{with} \quad \widehat R_\mu(t;h) =   R_\mu(t;h,\widehat H_t,\widehat L_t^2,  \widehat \sigma^2 (t) \color{black}),
\end{equation}
where $\widehat H_t,\widehat L_t^2$ are our local regularity estimators, and $\widehat \sigma^2 (t)$ is a suitable estimator of the errors variance.  We will show that, under mild conditions, $ \widehat R_\mu(t;h)  /  R_\mu(t;h) = 1 + o_{\mathbb P}(1)$, uniformly with respect to $h\in\mathcal H_N$. As a consequence, the rate of $h^*_\mu$ will coincide with that of the minimizer of  $R_\mu(t;h)$. For showing this, and deriving the convergence of $\widehat \mu_N^* (t)$, we require the following additional assumptions. For now, we focus on the independent design case. The common design case will be discussed in Section \ref{sec:comm_design}.

\begin{assumptionH}
\item \label{H:polynomelocaux} The estimator $\widehat{X}_n(t,h)$ is the Nadaraya-Watson  estimator with 
non-negative, symmetric and bounded kernel $K$, supported in $[-1,1]$. Moreover, $\inf_{|u|\leq 1}K(u) >0$.

\item \label{H:cd_grid_H} The bandwidth set $\Hh_N$ is a grid of points with at most $(N\lambda)^c$ points, for some $c>0$, such that $\max \Hh_N \rightarrow 0$ and
$N\lambda \min \Hh_N /\log(N\lambda)\rightarrow \infty$. Moreover, $\log(N)/\log^2(\lambda) \rightarrow 0$. \color{black}

\item\label{H:boundMn}  Constants $\underline{c},\overline{c}>0$ exist such that, for any $N$,  $ \underline{c} \leq M/\lambda = M / \EE(M) \leq \overline{c} $.

\item \label{H:densityT} The  density $g$ of the observation points $\Tni$ is Hölder continuous, and constants $\underline{c}_g,\overline{c}_g$ exist such that $0<\underline{c}_g\leq g(t) \leq \overline{c}_g$,   $\forall t\in I.$

\item\label{H:HL_conv} The estimators of $(H_t,L_t^2)\in(0,1)\times (0,\infty)$   admit concentration  bounds as in Theorem \ref{thm:regularity:H} and Theorem \ref{thm:regularity:L} with $\varphi = \left(\log\lambda\right)^{-2}$ and 
$\psi  =   \left(\log\lambda\right)^{-1}$, respectively. 

\end{assumptionH}
The condition on a kernel bounded from below (e.g., the uniform kernel) in \assrefH{H:polynomelocaux}, and the condition \assrefH{H:boundMn} can be relaxed at the cost of more involved technical arguments.   Regarding \assrefH{H:densityT}, if the design  density $g$ vanishes at   $t$, the pointwise convergence rate at $t$ of any nonparametric estimator would be degraded, and our assumption prevents this. We conjecture that by construction our risk bound \eqref{eq:mu:risk-boundV} adapts to low design, but we leave the study of this aspect for future work. Finally, since we necessarily have $h> (N\lambda)^{-1}$ for every $h\in \Hh_N $, and, by the last part of \assrefH{H:cd_grid_H}, $(N\lambda)^{-1/\log^2(\lambda)} \rightarrow 1$, the assumption \assrefH{H:HL_conv} guarantees that replacing the exponent $H_t$ by its estimate does not change the rate of the risk bound. For  $\widehat L_t^2$, which appears as a factor in the risk bound, a slower concentration rate is sufficient.

\begin{theorem}\label{thm:mu}
Let  $t\in I$, and assume that \assrefH{H:stationarity} to \assrefH{H:Lpmapprox},	  and  \assrefH{H:polynomelocaux} to  
\assrefH{H:HL_conv} hold true.
Then, we have $h_{\mu}^* = \Oo_\PP \{  (N \lambda)^{-1/(1 + 2H_t)} \},$
and the estimator $\widehat \mu_N^*(t) = \widehat \mu_N (t;h_{\mu}^*)$ defined in \eqref{eq:hatmu} and \eqref{eq:mu:risk-minimization} satisfies
\begin{equation*}
	\widehat \mu_N^*(t) - \mu(t) = \Oo_\PP\left\{ (N \lambda)^{-\frac{H_t}{1 + 2H_t}}  + N^{-1/2}\right\}.
\end{equation*}
\end{theorem}

The rate of the optimal bandwidth $h_{\mu}^*$ and the rate of convergence of $\widehat \mu_N^*(t)$ coincide with those obtained by  \cite{Golovkine2021} in the i.i.d. case. Our mean function estimator achieves the  minimax rate derived by  \cite{cai2011} for the mean function estimation. This convergence rate is slower than the parametric rate $ \Oo_\PP( N^{-1/2})$ in the \emph{sparse regime}  ($\lambda^{2H_t} \ll  N$), and achieves the parametric rate in the  \emph{dense regime}  ($\lambda^{2H_t} \gg N$). 
See \cite{zhang2016sparse} for the terminology.

We next derive the pointwise asymptotic distribution of our adaptive mean function estimation.
Usually, the rate of convergence in distribution for a nonparametric curve estimator is given by the power $-1/2$ of the effective sample, and the limit has the mean corrected by a bias term. In our context,  the effective sample size is expected to be  given by $N\lambda $ times the bandwidth. Meanwhile, the rate of convergence of the mean function estimator cannot be faster than the parametric rate $N^{-1/2}$ which corresponds to the ideal situation where all $N$ curves are observed without error at $t$. In the following we show that the effective sample size is given by $P_N(t;h_N)$ which, by construction, adaptively accounts for the two aspects. 

In a functional data context, the regularity of the mean function is necessarily equal to or larger than
that of the sample paths. As a consequence, the minimax optimal rates for the mean
function estimation are given by the sample path regularity, see \cite{cai2011}. Hence, from the minimax optimality perspective, the rate of convergence in distribution for a mean
function estimator has to depend on $H_t$. Assuming a higher regularity than the true one makes the convergence in distribution break down, as the bias term will tend to infinity. 

\begin{theorem}\label{thm:mean-clt}
Let $t\in I$ and assume that  \assrefH{H:stationarity} to 
\assrefH{H:erreur_oublie} and \assrefH{H:polynomelocaux} to \assrefH{H:HL_conv}
hold true.  Let  $h_{N}\in\Hh_N$, $N\geq 1$, such that 
\begin{equation}\label{h_as_norma}
	(N\lambda)^{1/(2H_t+1)}h_N  \rightarrow 0  . 
\end{equation}
Moreover, 
\begin{equation}\label{eq:sigma_clt}
	\frac{\sigma^2(t)}{P_N(t;h_N)}\sum_{n=1}^{N}\pi_n(t;h_N) \left\{\sum_{i=1}^{M_n}W_{n,i}^2(t;h_N)\right\}
	\overset{\PP}{\longrightarrow} \Sigma(t) \in [0,\infty),
\end{equation}
and 
\begin{equation}
	\operatorname{Var}_{M,T} \left( \frac{1}{\sqrt{P_N(t;h)} }\sum_{n=1}^N  \pi_n (t;h) \{X_n(t)- \mu(t)\} \right) \overset{\PP}{\longrightarrow}	\mathbb{S}_\mu (t)\in(0,\infty). 
\end{equation}
Then $		\sqrt{P_N(t;h_{N})}\left\{\widehat \mu_N(t;h_{N}) - \mu(t)\right\}\overset{d}{\longrightarrow} \mathcal{N}\left(0, \mathbb{S}_\mu (t) +\Sigma(t) \right).$
	\end{theorem}
	
	Condition \eqref{h_as_norma} makes the bias term negligible and thus avoids the usual mean correction used in the nonparametric regression. In the dense regime case, if in addition $\lambda h_N \rightarrow \infty$, then $\Sigma(t)=0$ and the $\widehat \mu_N(t;h_{N}) $ has the same asymptotic distribution as the infeasible empirical mean function estimator obtained with $X_i(t)$, $1\leq i\leq N$.  As expected, in the sparse regime case, the rate of convergence in distribution, given by $\EE[P_N(t;h_N)]^{-1/2}$, is slower than $N^{-1/2}$. Moreover, in this case $\Sigma(t)=\sigma^2(t)$ and $\mathbb{S}_\mu (t) = \operatorname{Var}(X(t))$. 

	\subsection{Adaptive autocovariance function estimator}
	
	The nonparametric estimation of the lag-$\ell$ autocovariance function with $\ell\geq 1$ seems less explored in the literature. \cite{kokoszka2017inference} consider the case of fully observed, error-free sample paths  and derive asymptotic results for the empirical autocovariance functions. \cite{zhong2023statistical} consider $MA(\infty)$ FTS observed with error over a fixed grid of design points, and use splines to estimate the sample paths. Their grid size corresponds to a dense regime and allows then to show that the empirical lag-$\ell$ autocovariance function constructed from the smoothed curves is asymptotically equivalent to the infeasible one obtained from the true sample paths. We propose here a nonparametric estimator of the lag-$\ell$ autocovariance function, $\ell \geq 1$, with independent or common design, in a sparse or dense regime, and which adapts to the regularity of the process generating the FTS.

	Let $s,t \in I$, and $\ell$ be an integer greater than or equal to 1. Let 
	\begin{equation}\label{PNl_def}
P_{N,\ell}(s,t;h) = \sum_{n=1}^{N-\ell} \pi_n(s;h)\pi_{n+\ell}(t;h),
\end{equation}
be the number of pairs  $  (X_{n},X_{n+\ell})$ with at least one pair  $(T_{n,i},T_{n+\ell,k})$ in the rectangle $[s-h, s+h]\times[t-h, t+h]$. 
The adaptive lag-$\ell$ cross-product kernel estimator of $\gamma_\ell(s,t) = \EE\left[X_{n}(s) X_{n+\ell}(t)\right] $ is $\widehat \gamma^*_{N,\ell}(s,t) = \widehat \gamma_{N,\ell}(s,t;h_\gamma^*) $, with 
\begin{equation}\label{eq:hat_crossproduct}
\widehat \gamma_{N,\ell}(s,t;h) \!=\!\! \sum_{n=1}^{N - \ell}\frac{\pi_n(s;h)\pi_{n+\ell}(t;h)}{P_{N,\ell}(s,t;h)}   \widehat X_n(s;h) \widehat X_{n+\ell}(t;h),
\end{equation}
where  $\widehat X_n(s;h) $ and $\widehat X_{n+\ell}(t;h)$ are NW estimators of $  X_n(s) $ and $ X_{n+\ell}(t)$, respectively. The selected bandwidth is a data-driven, optimal  bandwidth defined as
\begin{equation}\label{eq:gamma:risk-minimization}
h_{\gamma}^* \in \underset{h \in \Hh_N}{\argmin}\; \widehat R_\gamma(s,t;h),
\end{equation}
where $\widehat R_\gamma(s,t;h)$ is the estimate of 
\begin{align}\label{R_ga_theo}
R_\gamma(s,t;h)&=
3\nu^2_2\left(X_{1+\ell}(t)\right) L_s^2 h^{2H_s}\mathbb{B}(s|t;h,2H_s,0) + 3\nu^2_2\left(X_1(s)\right) L_t^2 h^{2H_t} \mathbb{B}(t|s;h,2H_t,\ell)\nonumber \\
&+ 3\left\{ \sigma^2(s) \nu^2_2(X_{1+\ell}(t)) 	 \mathbb{V}_{\gamma,0}(s,t;h) \color{black}+ \sigma^2(t)\nu^2_2(X_1(s))   \mathbb{V}_{\gamma,\ell}(s,t;h)   \right\} \nonumber\\
&+ 3 \sigma^2(s)\sigma^2(t) \mathbb{V}_\gamma(s,t;h) + \mathbb{D}(s,t;h)/P_{N,\ell}(s,t;h),
\end{align}
where for any $h > 0$, $\alpha > 0$, and any integer $\ell^\prime\geq 0$
\begin{align*}
\mathbb{B}(t|s;h,\alpha,\ell^\prime\!) &\!=\!\! \sum_{n=1}^{N-\ell}\!\frac{\pi_n(s;h)\pi_{n+\ell}(t;h)}{P_{N,\ell}(s,t;h)}  b_{n+\ell^\prime}(t;h,\alpha), \ \, b_n(t;h,\alpha) =\!\sum_{i=1}^{M_n}\!\left|\frac{\Tni - t}{h}\right|^{\alpha}\!\!\!\!  W_{n,i}(t;h),\\ 
\mathbb{V}_{\gamma,\ell^\prime}(s,t;h) &\!= \dfrac{1}{P_{N,\ell}(s,t;h)} \sum_{n=1}^{N-\ell} \frac{\pi_n(s;h)\pi_{n+\ell}(t;h)}{P_{N,\ell}(s,t;h)} 	 \!\max_{1\leq i \leq M_{n + \ell^\prime}}   W_{n+\ell^\prime ,i}(t;h) ,\\
\mathbb{V}_\gamma(s,t;h) &\!=\frac{1}{P_{N,\ell}(s,t;h)}\!\!\sum_{n=1}^{N - \ell}\! \frac{\pi_n(s;h)\pi_{n+\ell}(t;h)}{P_{N,\ell}(s,t;h)}\!  \max_{1\leq i \leq M_{n}}\!\! W_{n,i}(s;h)\!\! \max_{1\leq k \leq M_{n+\ell}}\!\! W_{n+\ell,k}(t;h),\\
\mathbb{D}(s,t; h) &\!=\! \EE(\!X_0\!\otimes\! X_{\ell}\! - \gamma_\ell\!)^2\!(\!s,\!t\!)	\!+\!2\!\!\!\!\!\!\sum_{k=1}^{N-\ell-1}\!\!\!\!p_k(s,t;h)\EE(\!X_0\!\otimes\! X_{\ell} - \gamma_\ell\!)(\!X_k\!\otimes\! X_{k+\ell}\! - \gamma_\ell\!)(\!s,t\!),\\
\text{where} &\quad p_k(s,t;h) = \sum_{i=1}^{N-k-\ell} \frac{\pi_i(s;h)\pi_{i+k}(s;h)\pi_{i+\ell}(t;h)\pi_{i + \ell + k}(t;h)}{P_{N,\ell}(s,t;h)}.
\end{align*}
Here, for any $f$ and $g$ real-valued functions,  $ 	(f\otimes g)(s,t) = f(s)g(t)$.

We will show that $	2R_\gamma(s,t;h)$ is a sharp bound of the quadratic risk $\EEMT\!\{\widehat \gamma_\ell(s,t;h) - \gamma_\ell(s,t)\}^2$, on a grid $\Hh_N$ of bandwidths. Like for the mean function estimation, the feasible bound $\widehat R_\gamma(s,t;h)$ is obtained by replacing $H$, $L^2$ and $\sigma^2$ values by the estimates introduced above. Moreover, the variances $\nu^2_2\left(X_{1}(s)\right)$ and $\nu^2_2\left(X_{1+\ell}(t)\right)$ are simply obtained as empirical variances of the presmoothing estimator $\{\widetilde X_n\}$ from Section \ref{sec:local-reg-estimator_nd}. Concerning $\mathbb{D}(s,t; h)$, let us first note  that under the $\LL^4_\Cc-m$-approximation assumption of the process $\{X_n, n \geq 1\}$, the autocovariances of the series $\{X_n\otimes X_{n + \ell}(s,t), n\geq 1 \}$ are absolutely summable. Similarly to mean estimation, this means that  we can simply take absolute values and bound $\mathbb{D}(s,t, h)$ by a constant. Details are provided in Lemma~\ref{lem:gamma:risk-bound}. 
For now, we assume that $\mathbb{D}(s,t; h)$ is given.  

We focus on the case of independent design, the case of common design is discussed in Section \ref{sec:comm_design}. To derive the asymptotic result for $h_\gamma^*$ and $\widehat \gamma^*_{N,\ell}(s,t)$, we add the following assumption on the bandwidth range.
\begin{assumptionH}
\item \label{H:cd_grid_H_autocov} $N(\lambda \min \Hh_N)^2 /\log(N\lambda) \longrightarrow \infty$.
\end{assumptionH}

\begin{theorem}\label{thm:gamma}
Assume the conditions \assrefH{H:stationarity} to \assrefH{D:smooth2}, \assrefH{H:Lpmapprox} for $p \geq 8$, \assrefH{H:polynomelocaux} to
\assrefH{H:densityT},
\assrefH{H:HL_conv} for $s,t\in I$,  and \assrefH{H:cd_grid_H_autocov} hold true.
Moreover, assume that a constant $	\mathfrak C >0$ exists such that 
\begin{equation}\label{mom4_mom2}
	\EE(X(u)-X(v))^4 \leq \mathfrak{C} \left[\EE(X(u)-X(v))^2 \right]^2,\qquad \forall u,v\in I.
\end{equation}
\color{black}
Let $H(s,t) = \min\{H_s, H_t\}$.
Then,  
\begin{equation*}
	h_\gamma^* = O_\PP \left(\max\left\{ (N\lambda^2)^{-\frac{1}{2\{H(s,t) + 1\}}}, (N\lambda)^{-\frac{1}{2H(s,t) + 1}} \right\}\right),
\end{equation*}
and 
\begin{equation*}
	\widehat \gamma^*_{N,\ell }(s,t) - \gamma_\ell(s,t) = 
	\Oo_\PP\left((N\lambda^2)^{-\frac{H(s,t)}{2\{H(s,t) + 1\}}} + (N\lambda)^{-\frac{H(s,t)}{2H(s,t) + 1}} + N^{-1/2}\right).
\end{equation*}
\end{theorem}

Let us note that
$$
(N\lambda^2)^{-\frac{1}{2\{H(s,t) + 1\}}} \gg  (N\lambda)^{-\frac{1}{2H(s,t) + 1}} \quad \Longleftrightarrow \quad 
\lambda ^{2H(s,t)} \ll  N,
$$ 
and 
$$
(N\lambda)^{-\frac{H(s,t)}{2H(s,t) + 1}} \ll  N^{-1/2} \quad \Longleftrightarrow \quad \lambda ^{2H(s,t)} \gg  N.
$$
As a consequence, in the `sparse' regime (\emph{i.e.}, $\lambda ^{2H(s,t)} \ll  N$), 
$$
\max \left\{ \left|\widehat\mu_N^*(s) -\mu(s)	 \right|, \left|\widehat\mu_N^*(t) -\mu(t)	 \right|  \right\} = o_\PP \left( \left|		\widehat \gamma^*_{N,\ell }(s,t) - \gamma_\ell(s,t)  \right|  \right).
$$
Moreover, the  convergence rates of  $\widehat \gamma^*_{N,\ell} (s,t)$,  $\widehat\mu_N^*(s)$ and $ \widehat\mu_N^* (t)$ are all slower than the parametric rate $O_\PP(N^{-1/2})$. In the `dense' regime  (\emph{i.e.}, $\lambda ^{2H(s,t)} \gg  N$), the three estimators attain the parametric rate.  As a  consequence,  the estimator $\widehat \Gamma^* _{N,\ell} (s,t) = \widehat \gamma^*_{N,\ell} (s,t) - \widehat\mu_N^*(s) \widehat\mu_N^* (t)$ of the autocovariance function estimator  $\Gamma_\ell (s,t) $ has the same  convergence rate as $\widehat \gamma^*_{N,\ell }(s,t) $.

The pointwise convergence rate for the lag$-\ell$ autocovariance function, obtained in Theorem \ref{thm:gamma}, coincides with the pointwise rate for the estimation of the covariance function, as obtained by \cite{Golovkine2021} in the i.i.d. case. This rate is given by the lowest regularity exponent $H$ at $s$ and $t$.

\subsection{The common design case}\label{sec:comm_design}

As noted by \cite{Golovkine2021}, the local bandwidth selection rules  defined in \eqref{eq:mu:risk-minimization} and \eqref{eq:gamma:risk-minimization} can be used for the  mean and covariance function estimation,  with both independent and common design. In the case of common design, where $\Tni\equiv T_i$, $1\leq i \leq \lambda$, the indicators $\pi_n(t;h)$  no longer depend on $n$, and they are all equal either to 0 or 1. That means that  $h_\mu^* $ and $	h_\gamma^* $ are automatically chosen in the set of admissible bandwidths where the $\pi_n(t;h)$ are all equal to 1. That also means that  the penalty terms  ${\mathbb{D}}_\mu(t;h)/P_N(t;h)$  and $\mathbb{D}(s,t;h)/P_{N,\ell}(s,t;h)$ can be removed from the risk bounds, because  they are constant on the range of admissible bandwidth values, and the risk bounds minimization is constrained to the admissible set. If the common design is equidistant, $h$ cannot be smaller than $1/\lambda$. For both mean and autocovariance functions, two cases can occur: the minimum of the risk bound without the penalty term is attained in the interior of the admissible set of $h$ (dense regime case), or on the left boundary where the bias term will be larger than the variance term (sparse regime case). Thus, our kernel smoothing automatically selects between linear interpolation and smoothing by choosing the optimal bandwidth in a data-driven manner. This  is illustrated in our real data analysis for the mean function estimation. As a consequence of these facts, we can deduce the following result for which the justification is obvious and is thus omitted.

\begin{theorem}\label{thm:com_des}
Assume that  $\Tni$ belong to a common design as in condition (H\ref{H:Tni}).
\begin{enumerate} 
	\item Assume that the conditions of Theorem \ref{thm:mu} are satisfied and
	$\widehat \mu_N^*(t) = \widehat \mu_N (t;h_{\mu}^*)$ with $h_\mu^*$ defined as in \eqref{eq:mu:risk-minimization}. Then, 
	$
	\widehat \mu_N^*(t) - \mu(t) = \Oo_\PP( \lambda^{-H_t} + N^{-1/2}).
	$

	\item Assume that the conditions of Theorem	\ref{thm:gamma} are satisfied and $	\widehat \gamma^*_{N,\ell}(s,t) = \widehat \gamma_{N,\ell}(s,t;h_\gamma^*)  $ with $h_\gamma^*$ defined as in \eqref{eq:gamma:risk-minimization}.
	Then,  
	$
	\widehat \gamma^*_{N,\ell }(s,t) - \gamma_\ell(s,t) = 
	\Oo_\PP(\lambda^{-H(s,t)} + N^{-1/2}).
	$
\end{enumerate}

\end{theorem}


\section{Numerical study}\label{sec:emp_study}
This section presents a  Monte Carlo study and an application to daily voltage curves from the Individual Household Electricity Consumption dataset  \citep[]{misc_individual_household_electric_power_consumption_235}. 
The  results were obtained using an R package 
which is publicly available at \url{https://github.com/hmaissoro}. The Epanechnikov kernel ($K(u)=(3/4)(1-u^2)$ for $|u|\leq 1$, and $0$ otherwise) was used in all experiments.

\subsection{Simulation setting}\label{sec:simu_exp}
We consider three types of FTS $\{X_n\}$ and investigate the effectiveness of our methods in the case of non-differentiable sample paths with different local regularity exponents. 
The three types of FTS that are considered are versions of a FAR(1) process with an associated innovation process $\{\xi_n\}$,
\begin{equation}\label{FTS_mod1}
X_n(u) = \mu(u) + \int_0^1\psi(u,s)(X_{n-1}(s)-\mu(s))ds + L_t\xi_n(u),
\end{equation}
where $\mu(t)= 4\sin(3\pi t/2)$, $\psi(u,s)=\kappa\exp(-(u+2s)^2)$ and the constant $\kappa$ is chosen so that the operator norm $\vertiii{\cdot}_\infty$ of the integral operator defined by $\psi(u,s)$ is equal to $0.5$, while $L_t=2$. A series of 100 burn-in steps are used to initialize \eqref{FTS_mod1}.

\begin{figure}[th]
\centering
\begin{tabular}{ccc}
	\includegraphics[width=0.3\textwidth]{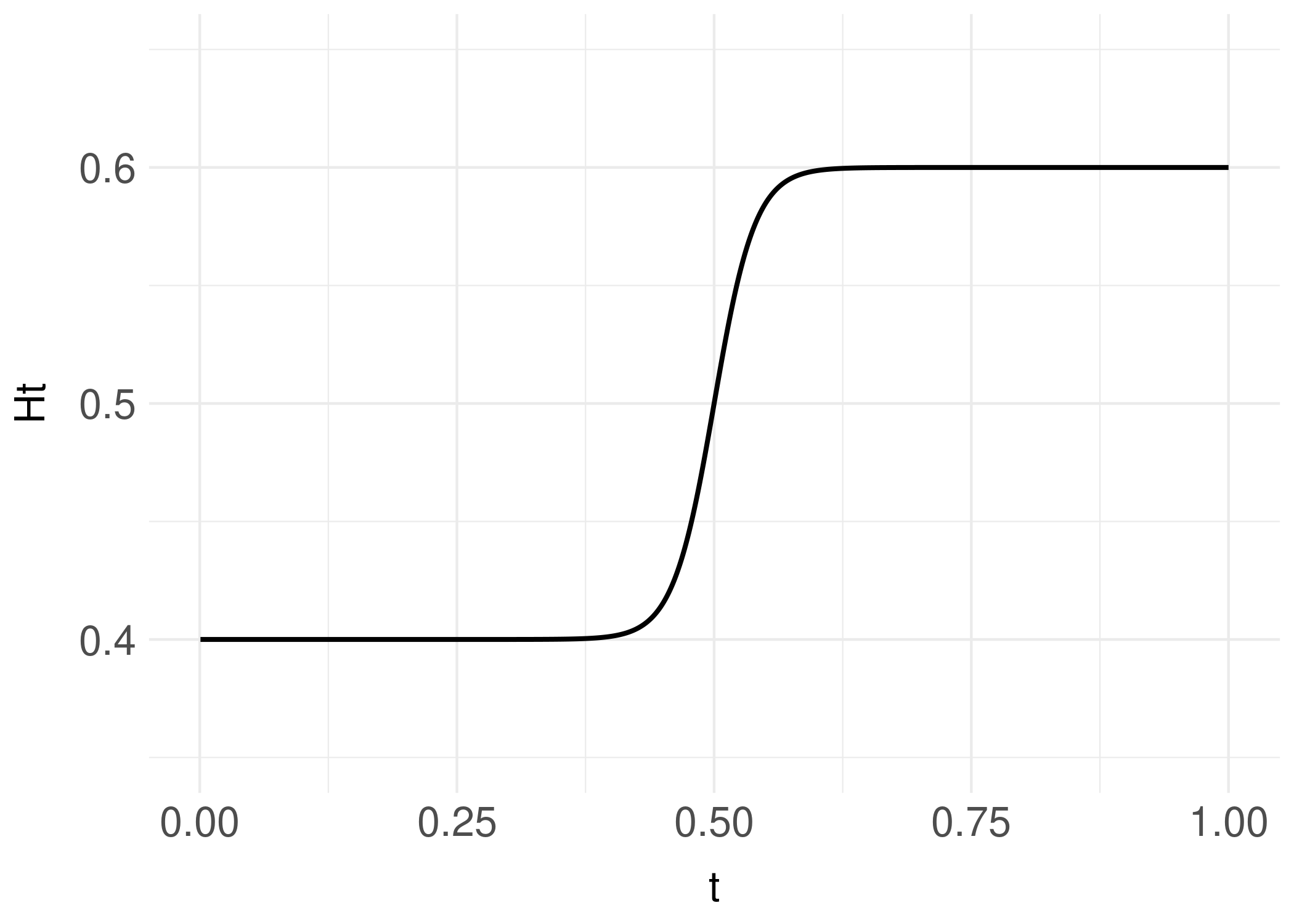} &
	\includegraphics[width=0.3\textwidth]{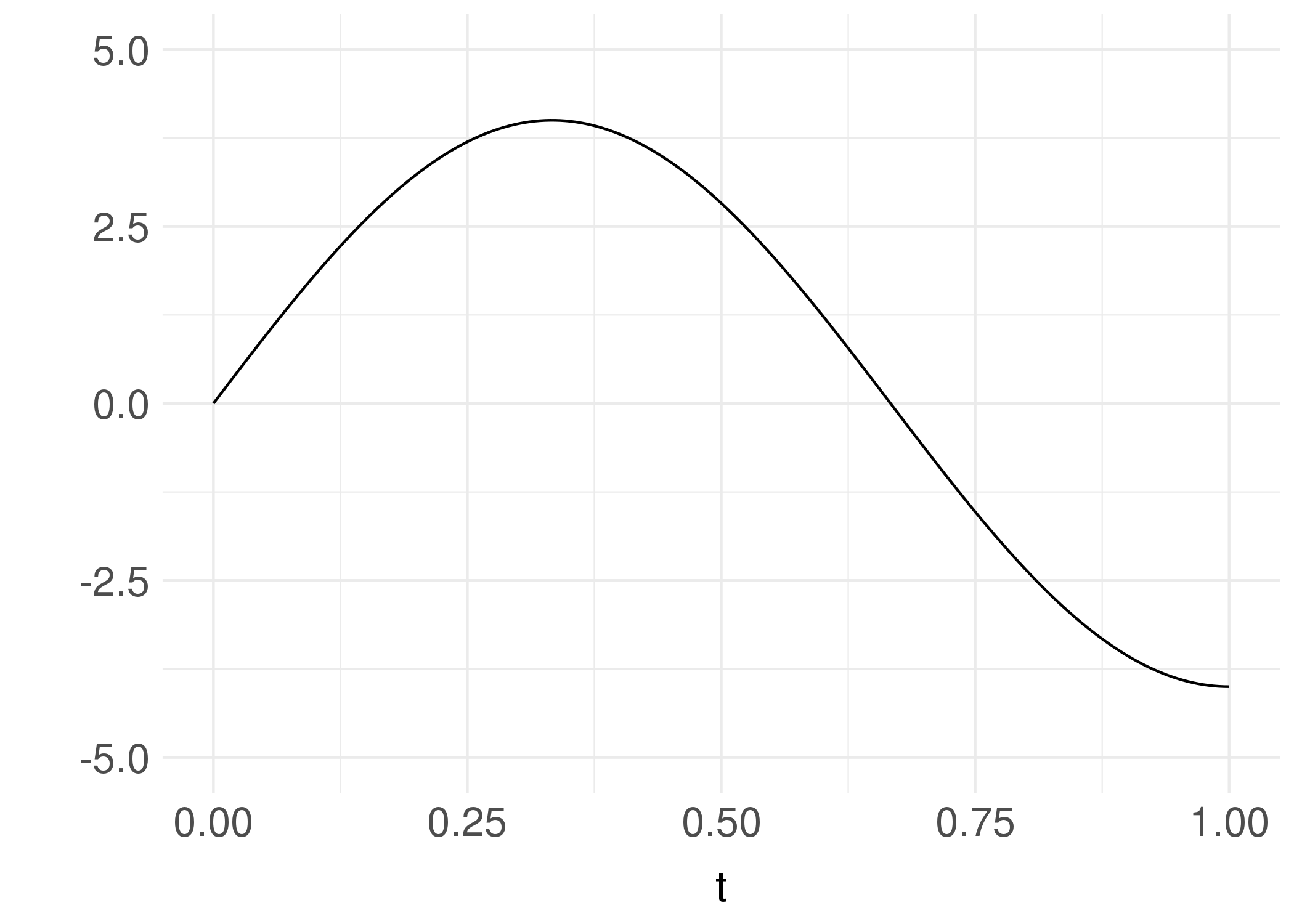} &
	\includegraphics[width=0.3\textwidth]{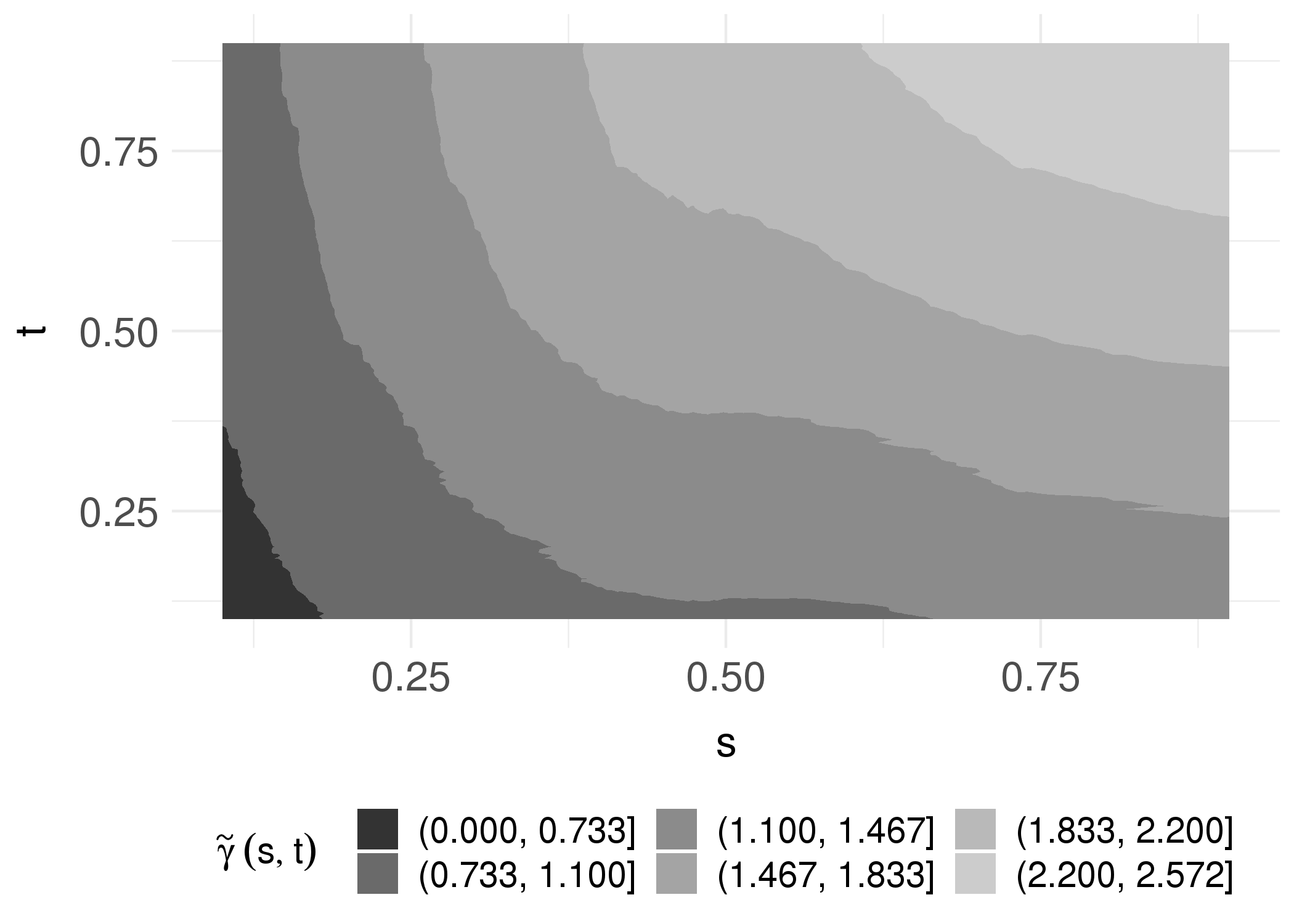}
\end{tabular}
\caption{\small Simulation parameters. \textbf{Left:} Logistic local exponent function $H_t$ used in FTS Model 2 and 3. \textbf{Middle:} The mean function $\mu$ used in FTS Model 1 and 2. 
	\textbf{Right:} The  empirical approximation of the lag$-1$ autocovariance function $\gamma_1(s,t)$ obtained  from a large sample in FTS Model 2 when $\mu\equiv 0$. }
\label{fig:true_locreg_param_mean_kernel}
\end{figure}

The simulation results we present here are obtained with the series $\{X_n\}$ generated in what we call \textbf{FTS Model 2}~:  $\{X_n\}$ is a FAR(1) as in \eqref{FTS_mod1}, with $\{\xi_n\}$ independently generated from a MfBm with a logistic Hurst index function (see Figure~\ref{fig:true_locreg_param_mean_kernel}). \textbf{FTS Model 1} is a version of \textbf{FTS Model 2} with a constant Hurst index $H_t$ instead of the logistic one, while \textbf{FTS Model 3} is another version of \textbf{FTS Model 2} with the mean function $\mu(t)$ and the function $\psi(u,s)$ learned from the daily voltage of 
the Individual Household Electricity Consumption dataset \citep{misc_individual_household_electric_power_consumption_235}. The results obtained with \textbf{FTS Model 1} and \textbf{3}, as well as details of the setups of these models, are presented in \citet{hassan2024supp}.

%
%
%
%

To obtain the data points according to \eqref{eq:data-model}, 
the integers $M_n$ are randomly generated  uniformly between $0.8\lambda $ and  $1.2\lambda$, while the 
$\{T_{n,i}\}$ are uniformly distributed over $(0,1]$. The errors $\varepsilon_{n,i}$ are Gaussian with constant variance $\sigma^2 = 0.25^2$.
We consider $(N, \lambda)\in\{(150,40), (1000, 40), (400,300), (1000,1000)\}$.
For each 
setup, 
we generate $R=400$ independent  
series.

\subsection{Local regularity estimation} \label{sec:emp_study:loc_reg}

\begin{figure}[th]
\centering
\includegraphics[height=0.33\textheight,width=0.9\textwidth]{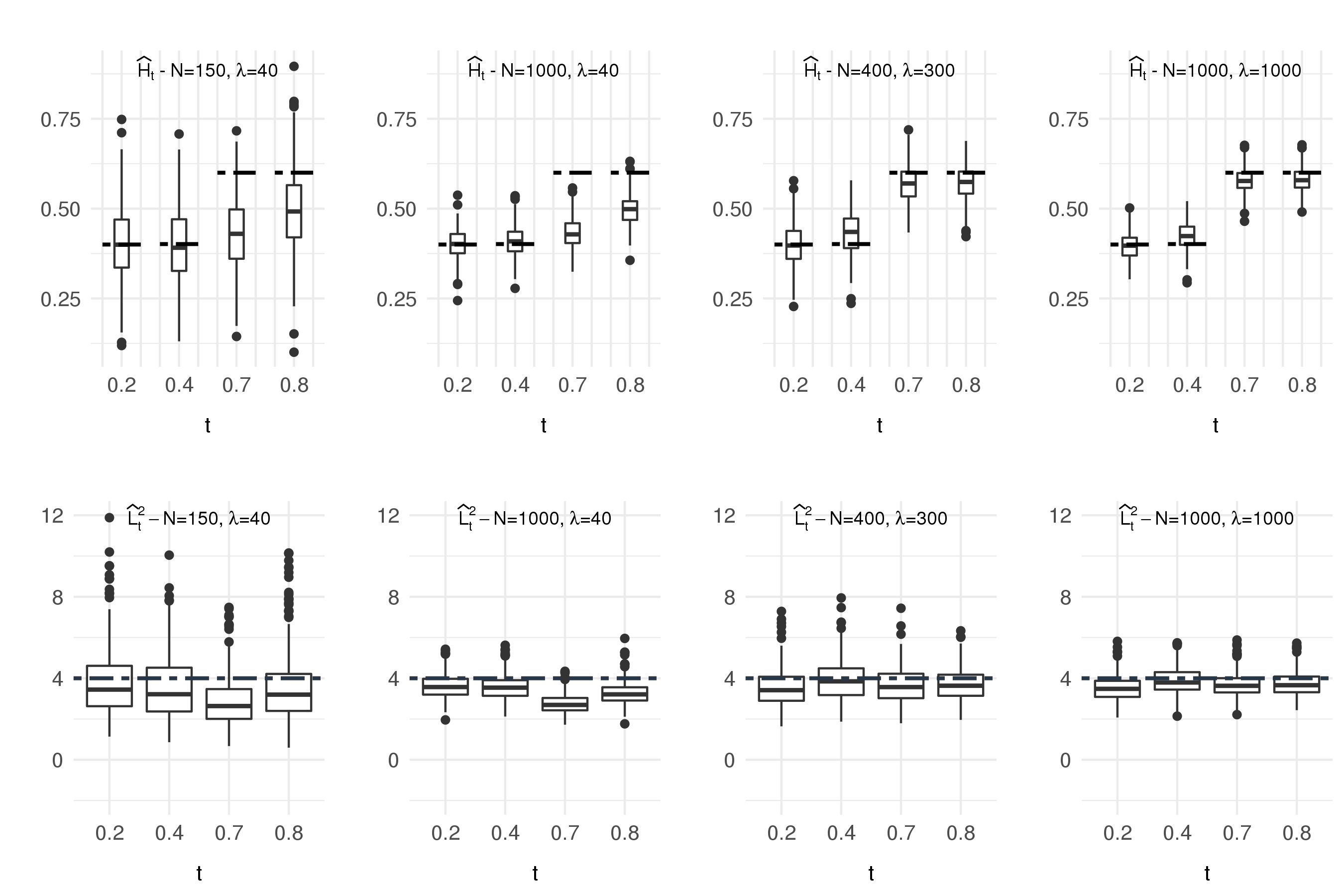}
\caption{\footnotesize Boxplots of $R=400$ pointwise estimates of $\widehat{H}_t$ and $\widehat{L}^2_t$, for $t\in\{0.2,0.4,0.7,0.8\}$ and four pairs $(N,\lambda)$, in  FTS Model 2. The dashed horizontal lines indicate the true values of $H_t$ and $L^2_t$.}
\label{fig:locreg_far_mfBm_d2}
\end{figure}

Our  approach for the estimation of $H_t$ and $L^2 _t$ depends on two tuning parameters: the window length $\Delta$ used in \eqref{eq:proxyH}, and the presmoothing bandwidth used in \eqref{eq:pre-smooth}. The presmoothing bandwidth is selected by a cross-validation procedure described in \citet{hassan2024supp}. Concerning $\Delta$, Theorems \ref{thm:regularity:H} and \ref{thm:regularity:L} propose the choice $\Delta = \exp\left( -(\log\lambda)^\gamma\right)$ for some $\gamma\in(0,1)$. On the basis of extensive simulations, for which the details are provided in \citet{hassan2024supp}, we set $\gamma = 1/3$.
Figure~\ref{fig:locreg_far_mfBm_d2} shows the boxplots of $\widehat{H}_t$ and $\widehat{L}_t^2$ defined in \eqref{eq:H0hat} for the four pairs $(N, \lambda)$ at four points  $t\in I=(0,1].$ The bias of the regularity parameters estimates  decreases as $\lambda$ increases, and the boxplot are more concentrated as $N$ increases. Overall, the local regularity estimators show good finite sample performance.

\subsection{Mean function estimation} 

\begin{figure}[th]
\vspace{-.1cm}
\centering
\includegraphics[width=0.85\textwidth, height=.41\textwidth]{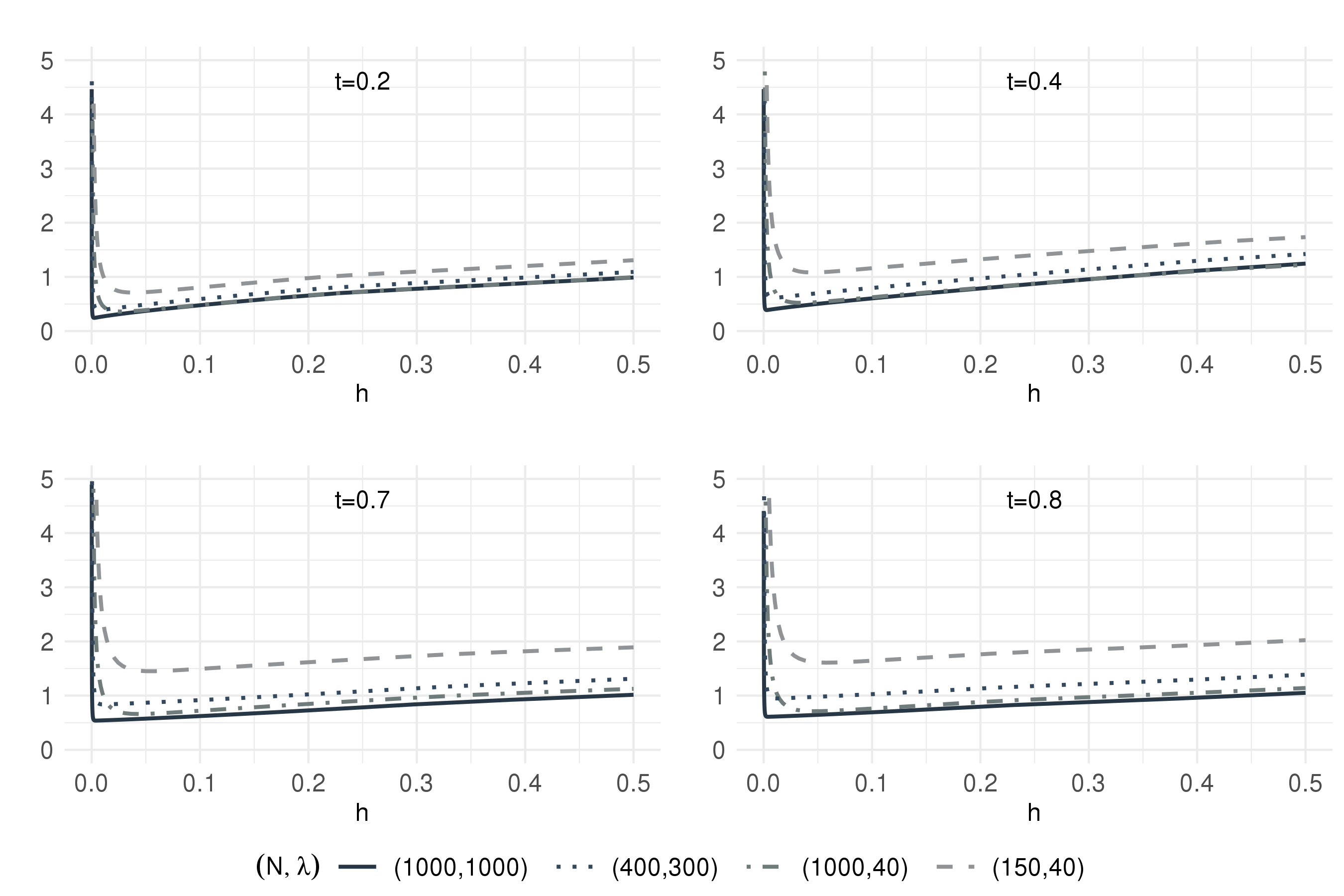}
\caption{\footnotesize Empirical average of the risk function $\widehat{R}_\mu(t;h)$  at $t \in\{0.2, 0.4, 0.7, 0.8\}$ over $400$ independent functional time series generated according to FTS Model 2, with four setups $(N,\lambda)$.}
\label{fig:mu_risk_far_mfBm_d2}
\end{figure}

Our adaptive `smooth first, then estimate' estimator of the mean function is constructed with the bandwidth $h_{\mu}^*$ defined as in \eqref{eq:mu:risk-minimization}, obtained by minimizing the estimated bound $2\widehat{R}_\mu(t;h)$ of the pointwise quadratic risk. Instead of the dependence coefficient $\mathbb{D}_\mu(t;h)$, we simply consider 
\begin{align*}
\overline{\mathbb{D}}_\mu (t;h) 
=& \frac{1}{N} \sum_{n=1}^{N} \!\left\{\!\widetilde{X}_n(t) \!- \widehat{\nu}_1(X(t))\!\right\}^2 
\\ &+ 2\sum_{\ell = 1}^{N-1}\frac{1}{N \!- \!\ell} \left| \sum_{n = 1}^{N - \ell - 1}\!\left\{\!\widetilde{X}_n(t) - \widehat{\nu}_1(X(t))\right\} \!\left\{\!\widetilde{X}_{n+\ell}(t) \!- \widehat{\nu}_1(X(t))\!\right\}\right|,
\end{align*}
with $\{\widetilde{X}_n\}$ the presmoothed curves as defined in \eqref{eq:pre-smooth} and $\widehat{\nu}_1(X(t))$ their empirical mean at $t$. Figure \ref{fig:mu_risk_far_mfBm_d2} presents the  average of the risk function $\widehat{R}_\mu(t;h)$ over $400$ independent time series generated according to FTS Model 2,  with four setups $(N,\lambda)$. The plots provide evidence that $h\to R_\mu(t;h)$ is a convex function which converges to zero as $N$ and $\lambda$ become larger.

\begin{figure}[th]
\vspace{-.25cm}
\centering
\includegraphics[width=0.8\textwidth, height=0.25\textheight]{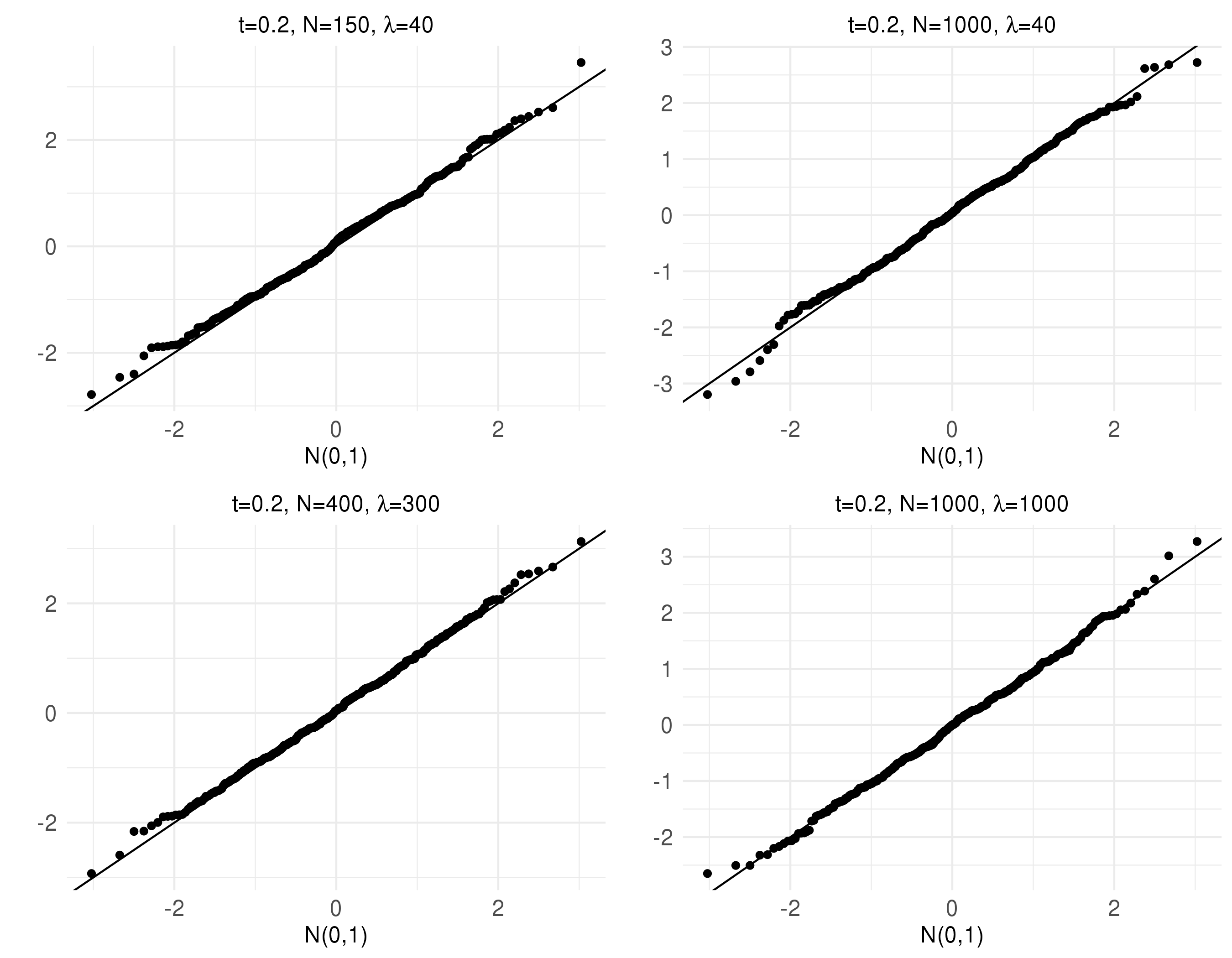}
\caption{\small Normal $Q-Q$ plots of $\sqrt{P_N(t;h_N)}\left(\widehat{\mu}_N(t;h_N) -\mu(t)\right) / \sqrt{\widehat{\mathbb{S}}_\mu (t) +\widehat{\Sigma}(t)}$ at  $t=0.2$, with $h_N = \{h_\mu^*\}^{1.1}$. Results obtained with $400$ independent  time series generated in the  FTS Model 2.}
\label{fig:mu_tcl_d2}
\end{figure}

Table~\ref{tab:mean_far_mfBm_d2} shows the bias and standard deviation of the estimates of $\widehat \mu_N^* (t) = \widehat \mu_N (t;h_{\mu}^*)$ obtained for functional time series generated according to the FTS Model 2. As expected the bias and the variance decrease as $N,\lambda \to \infty$. The estimated standard deviations increase as $t$ increases, which may be surprising given that the sample paths become smoother to the right of the domain $I$. However, larger $t$ also means larger $\operatorname{Var}(X_t)$ \citep[see][for the variance plot]{hassan2024supp}, and the consequence is less precise estimates of the mean. 
Finally, we study the asymptotic distribution of $\widehat \mu_N^* (t)$. The $Q-Q$ plots in Figure \ref{fig:mu_tcl_d2} show that the Gaussian limit, as stated in Theorem \ref{thm:mean-clt}, is an accurate approximation. Indeed, we notice that the distribution of $ P_N(t;h_N)^{1/2} \{\widehat \Sigma(t) + \widehat {\mathbb S}_\mu (t) \}^{-1/2}   \left\{\widehat{\mu}_N(t;h_N) -\mu(t)\right\}$ is close to the standard normal distribution for all $(N,\lambda)$ considered. The estimates $\widehat \Sigma(t)$ and $\widehat {\mathbb{S}}_\mu (t)$ are defined in \citet{hassan2024supp}.

\begin{table}[th]
\begin{center}
	{\footnotesize
		\begin{tabular}{|cc|cc|cc|cc|cc|}
			\hline
			&&\multicolumn{2}{|c|}{\bfseries $t = 0.2$}&\multicolumn{2}{|c|}{\bfseries $t = 0.4$}&\multicolumn{2}{|c|}{\bfseries $t = 0.7$}&\multicolumn{2}{|c|}{\bfseries $t = 0.8$}\\
			$N$&$\lambda$& Bias & Sd & Bias & Sd & Bias & Sd & Bias & Sd \tabularnewline
			\hline
			$ 150$&$  40$ &0.0056 &0.2079 &0.0112&0.2692 &0.0329&0.3259  &0.0497&0.3417\\
			$1000$&$  40$ &0.0005 &0.0883 &-0.0062&0.1139&0.0119&0.1353  &0.0213&0.1425\\
			$ 400$&$ 300$ &0.0074 &0.1283 &0.0049&0.1626 &0.0119&0.1944  &0.0150&0.2044\\
			$1000$&$1000$ &-0.0020&0.0849 &0.0004&0.1094 &-0.0003&0.1301 &0.0003&0.1369\\
			\hline
	\end{tabular} }
\end{center}
\caption{Bias and standard deviation (Sd) of the mean function estimates obtained from $400$ independent  time series generated in the FTS Model 2.} \label{tab:mean_far_mfBm_d2}
\end{table}

\begin{figure}[th]
\centering
\includegraphics[height=0.31\textheight,width=0.9\textwidth]{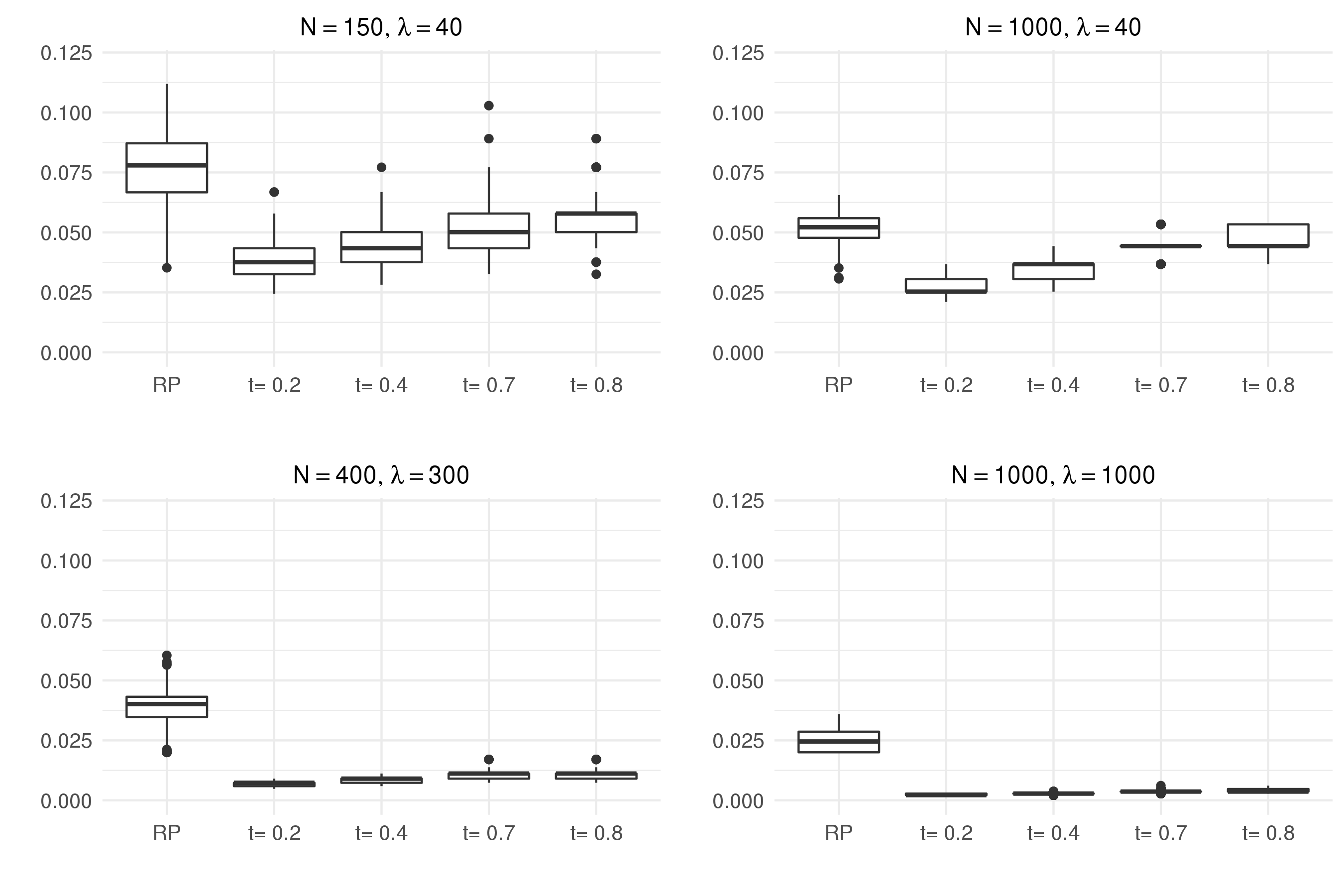}
\caption{\small Bandwidths selected by \texttt{RP20} (left boxplot)  and by our   local approach   for the mean  estimation at $t\in\{0.2,0.4,0.7,0.8\}$; results  from $400$ independent  series generated in the FTS Model 2.}
\label{fig:mu_bandwidth_d2}
\end{figure}

We conclude this section  with a comparison with the  procedure of \citet{rubin2020sparsely}, procedure  referred to as \texttt{RP20}, in the context of the FTS Model 2. A similar comparison in the context of the FTS Model 3 can be found in \citet{hassan2024supp}. 
\citet{rubin2020sparsely} proposed a locally linear estimator of the mean function, which we  denote by $\widehat{\mu}_{\texttt{RP}}$, in sparsely observed settings. Their bandwidth is selected by $K$-fold cross-validation using the Bayesian optimisation algorithm implemented in MATLAB.
The implemented procedure is such that the observations times $\{T_{ij}\}$ are randomly sampled over a regular discrete grid of $241$ points. In addition, since the implementation of $K$-fold cross-validation is time consuming, a projection on a B-spline basis is proposed for dimension reduction in the Bayesian optimisation. In Figure \ref{fig:mu_bandwidth_d2} we present the boxplots of the selected bandwidths according to \texttt{RP20}'s global approach and to our local approach. The selected bandwidths have comparable sizes in almost all setups $(N,\lambda)$. As expected from the increasing shape of the function $H$, our local bandwidths are smaller for $t$ in the first half of $I$ and increase as $t$ gets closer to 1. Table \ref{tab:mean_mse_far_mfBm_d2} shows the ratio of the Monte-Carlo estimates of the Mean Square Error (MSE) of our mean function  estimator and the \texttt{RP20}'s local linear estimator. Although the ratio is close to 1, our estimator shows slightly better performance (ratio less than 1) in almost all setups.


\begin{table}[th]
\vspace{-.4cm}
\begin{center}
	{\footnotesize \begin{tabular}{|cc|c|c|c|c|}
			\hline
			{$N$}&{$\lambda$}& {\bfseries $t = 0.2$}&{\bfseries $t = 0.4$}& {\bfseries $t = 0.7$}&{\bfseries $t = 0.8$}\\
			\hline
			$ 150$&$40$   & 0.9689 & 0.9321 & 0.9520 & 0.9537 \\
			$1000$&$40$   & 0.9710 & 0.9414 & 0.9228 & 0.9208 \\
			$ 400$&$ 300$ & 1.0131 & 0.9716 & 0.9959 & 0.9867 \\
			$ 1000$&$1000$& 0.9914 & 1.0015 & 0.9917 & 0.9949 \\
			\hline
	\end{tabular} }
	\caption{\small MSE ratio for our mean estimator and \texttt{RP20}; results from $400$   series generated in FTS Model 2.}
	\label{tab:mean_mse_far_mfBm_d2}
\end{center}
\end{table}

\vspace{-.4cm}
\subsection{Autocovariance function estimation}
To focus on the specific aspects related to the estimation of the lag-$\ell$ autocovariance function,  we consider series generated as in  FTS Model 2 but with the mean function set equal  to zero. We set $\ell = 1$. 
An accurate approximation of $\gamma_1(s,t) = \EE[X^{(n)}_sX^{(n-1)}_t]$ is shown in Figure~\ref{fig:true_locreg_param_mean_kernel} \cite[see][for details]{hassan2024supp}.
In this case, $\widehat \Gamma^* _{N,1} (s,t) = \widehat \gamma_{N,1} (s,t;h_\gamma^*)$, 
with $\widehat \gamma_{N,1} (s,t;h)$ 
defined in  \eqref{eq:hat_crossproduct} and the bandwidth $h_\gamma^*$ obtained from \eqref{eq:gamma:risk-minimization}. Further details on the optimization in  \eqref{eq:gamma:risk-minimization} are given in \citet{hassan2024supp}. The results of the estimation of the lag$-1$ autocovariance function for two setups $(N,\lambda)$  are presented in the Table~\ref{tab:autocov_far_mfBm_d4}. Larger Sd values occur for smaller $N$ and/or for points $(s,t)$ with larger values of $\gamma_1(s,t)$.

\begin{table}[th]
\begin{center}
	{\footnotesize	\begin{tabular}{|cc|cc|cc|cc|cc|}
			\hline
			&&\multicolumn{2}{|c|}{\bfseries $(s,t) = (0.2,0.4)$}&\multicolumn{2}{|c|}{\bfseries $(s,t) = (0.4,0.7)$}&\multicolumn{2}{|c|}{\bfseries $(s,t) = (0.7,0.8)$}&\multicolumn{2}{|c|}{\bfseries $(s,t) = (0.8,0.2)$}\\
			$N$&$\lambda$& Bias & Sd & Bias & Sd & Bias & Sd & Bias & Sd \tabularnewline
			\hline
			$ 150$&$40$& 0.0019&0.3359 &0.0307&0.5193  &0.0371&0.6675 &0.0102&0.4058\\
			$1000$&$40$& 0.0052&0.1303 &-0.0004&0.1893 &0.0026&0.2398 &0.0126&0.1568\\
			\hline
	\end{tabular}  }
	\caption{\small Bias and standard deviation (Sd) of the lag-$1$ cross-product function $\gamma_1(s,t)$ estimation in FTS Model 2 when $\mu\equiv 0$; results obtained from $400$ independent series.}
	\label{tab:autocov_far_mfBm_d4}
\end{center}
\end{table}

\vspace{-.5cm}
\subsection{Real data analysis}

Predicting electrical energy consumption is essential for planning electricity production and significantly reduces the problems of storage and overproduction. A key step in this objective is to be able to accurately estimate  the evolution  of the electricity production parameters (such as the voltage), and functional time series are an effective approach for this purpose. To illustrate, we consider the data provided by the Individual Household Electricity Consumption dataset from the UC Irvine Machine Learning Repository
\citep{misc_individual_household_electric_power_consumption_235}. 
It contains various measurements of electricity consumption in a household near Paris, with a sampling rate of one minute from December 2006 to November 2010. The data of interest here are $1358$ voltage curves with a common design of $1440$ points (corresponding to minute-by-minute observations), normalized so that $I= (0, 1]$.  There are about $5.8\%$ daily curves missing from the dataset, but we decided to neglect the missingness effect and consider the series as complete.

\begin{figure}[ht]
\centering
\begin{tabular}{cccc}
	\small (a) $\widehat{\mu}_N^*$ & \small (b) $\widehat{\mu}_{\texttt{RP}}$\\
	\includegraphics[height= 0.15\textheight,width=0.47\textwidth]{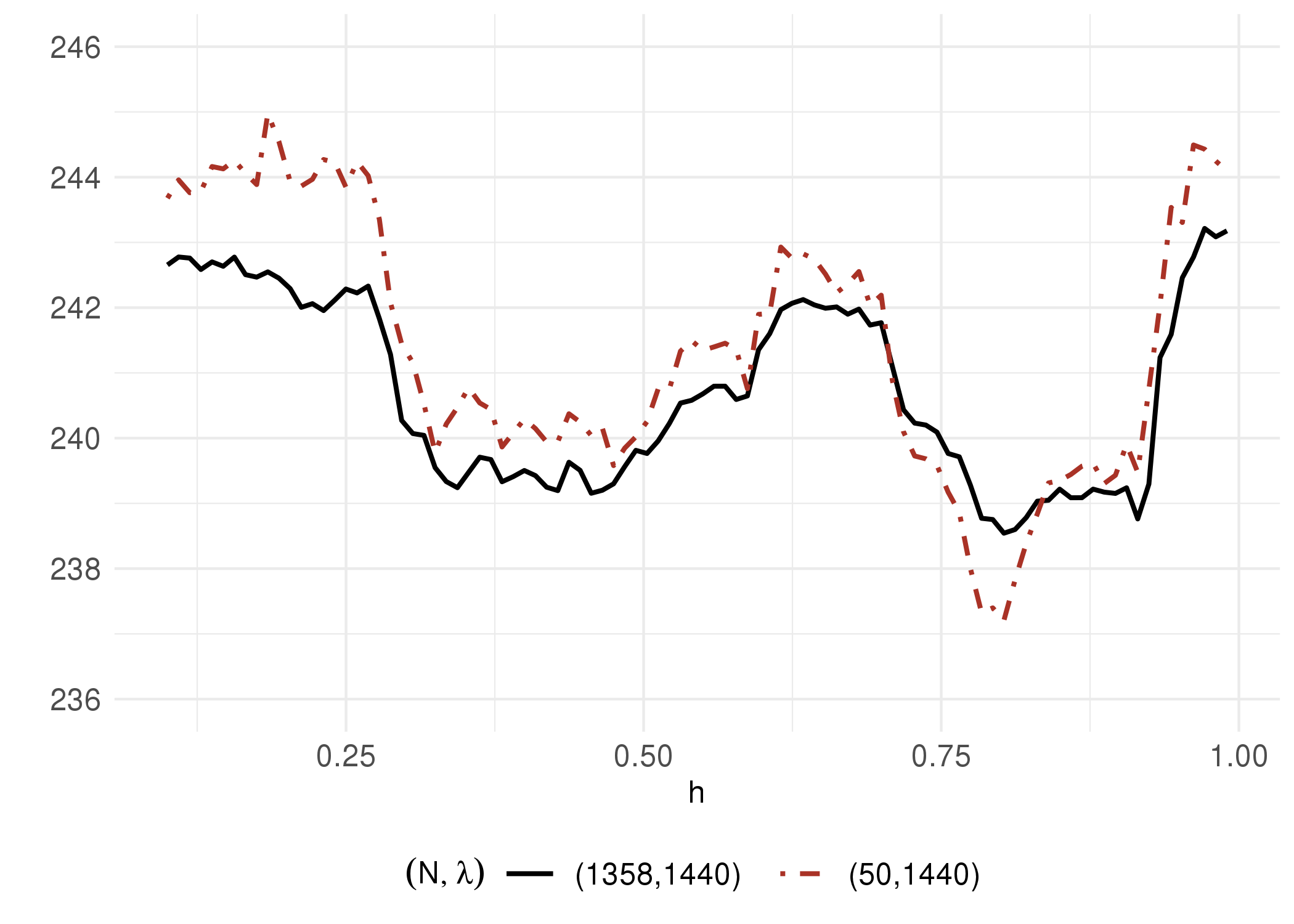} &
	\includegraphics[height= 0.15\textheight,width=0.47\textwidth]{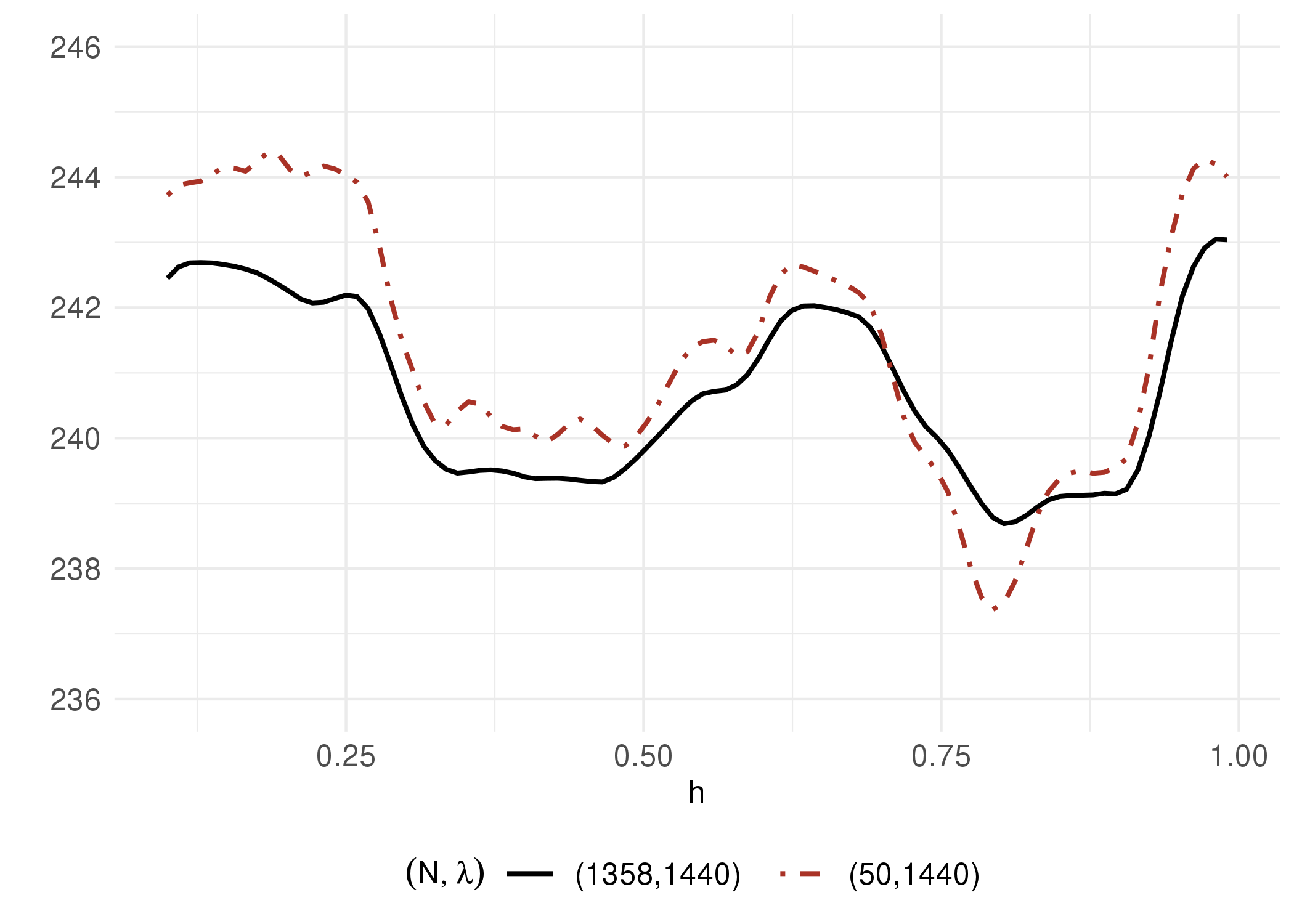}\\
	\small (c) $t\to \widehat{H}_t$ & \small (d) $t\to \widehat{L}_t^2$\\
	\includegraphics[height= 0.15\textheight,width=0.47\textwidth]{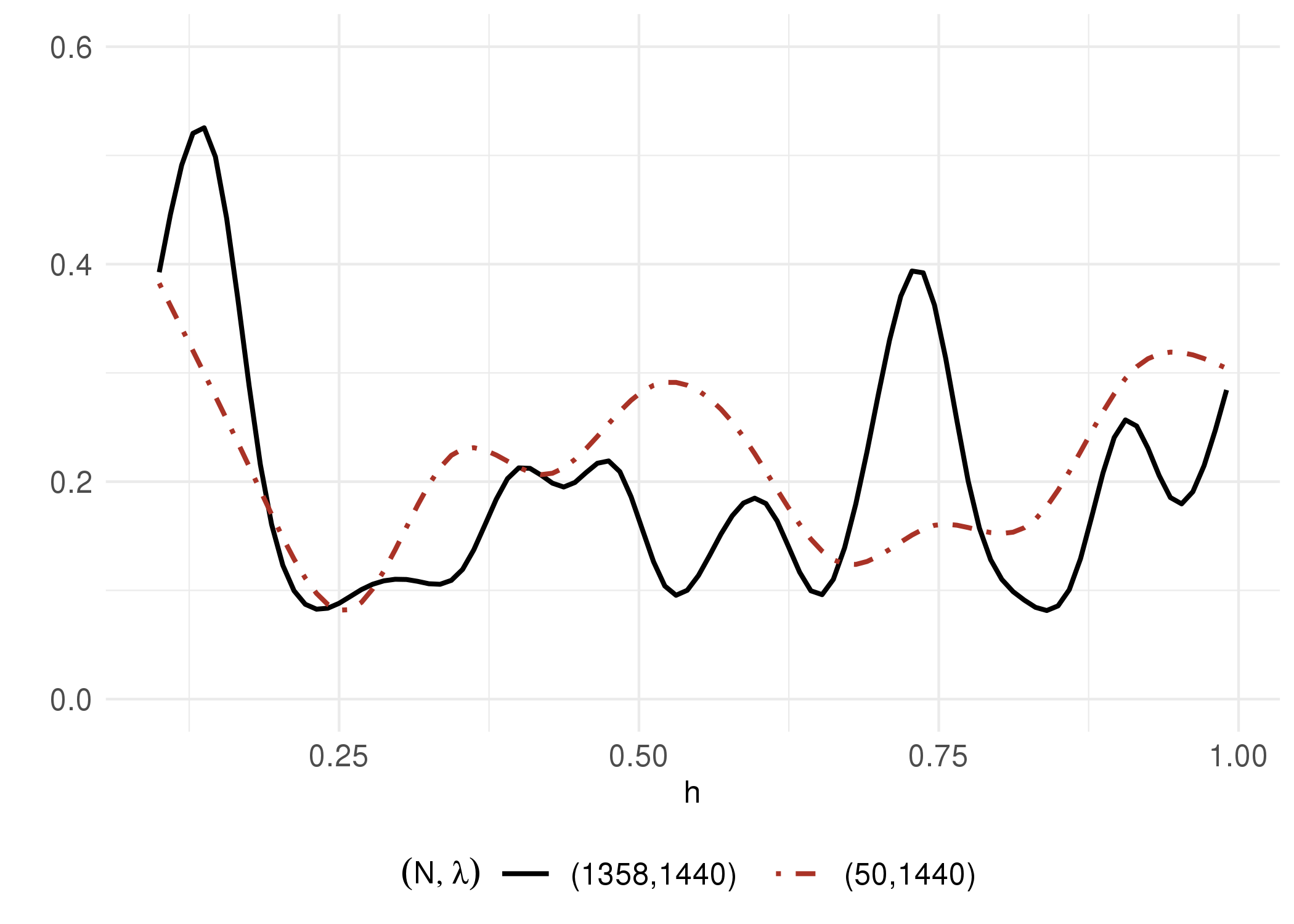} &
	\includegraphics[height= 0.15\textheight,width=0.47\textwidth]{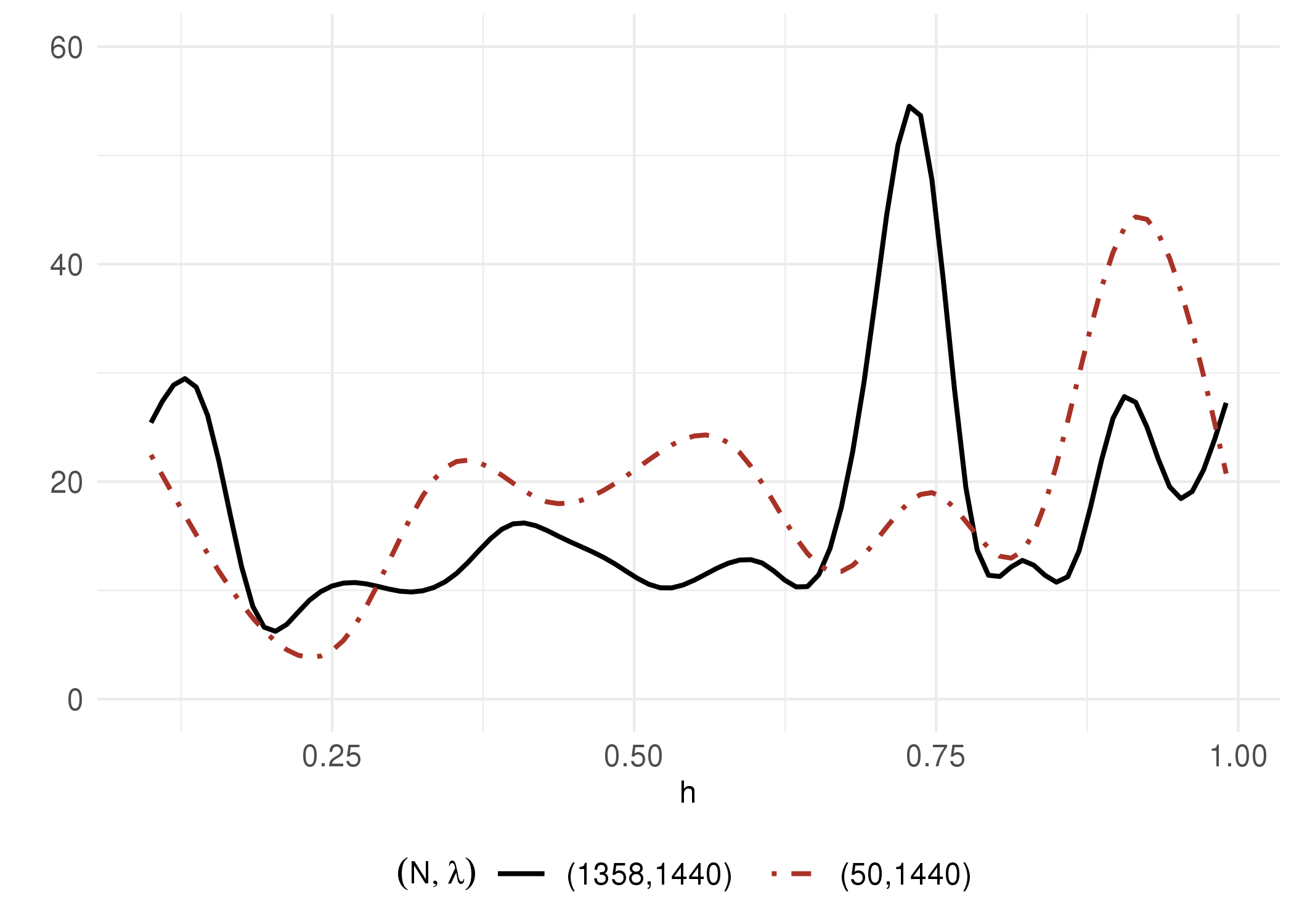}
\end{tabular}
\caption{\footnotesize Estimation of the mean function with our procedure and with \texttt{RP20} for the daily voltage curves, $N=1358$ (solid black lines) and $N=50$ (dotted red lines)~: (a) the estimation of $\mu$ with our procedure; (b) the estimation of $\mu$ with \texttt{RP20}; (c) and (d) the estimations of the local regularity parameters.  }
\label{fig:daily_voltage}
\end{figure}

Figure \ref{fig:daily_voltage} shows the estimates of the daily mean voltage curve with our procedure and the procedure \texttt{RP20}, on the grid of $1440$  points. We have considered the full series ($N=1358$) and a sub-series of $N=50$ consecutive curves from a period with no missing days. As  expected, our estimate is more irregular than that obtained by the \citet{rubin2020sparsely} procedure. It is worth noting that when $N=1358$, for most of the $1440$ points over the fixed grid, our optimal bandwidth leads to degenerate smoothing using only one data point, which is equivalent to interpolation. This is no longer the case when $N=50$, where the larger bandwidths lead to smoothing using up to 18 data points. In other words, our adaptive procedure automatically chooses between the interpolation and smoothing in the common design setting. In the common design with sparsely sampled curves (\emph{i.e.}, when $\lambda^{2H_t}\ll N$) interpolation is minimax rate optimal \citep[see][]{cai2011}. In our application the number of design points per curve and the full time series length are comparable (1440 versus 1358), and $H_t$ is much less than 1/2 for almost all $t$,  suggesting a sparse regime that our adaptive mean estimation approach automatically detects. With $N=50$, the setup is one of densely sampled curves, and our bandwidth is automatically chosen accordingly.


\section{Conclusions}\label{sec:conclusions}

We have studied a notion of local regularity for the process generating the sample paths of a stationary, weakly dependent functional time series (FTS). The paths are observed with heteroscedastic errors on discrete sets of design points, which may be fixed or random. The weak dependence condition we consider is satisfied by a large panel of FTS. Using a Nagaev-type inequality, we derive bounds on the concentration of the  regularity estimators. Using the regularity estimators,  adaptive mean and autocovariance function nonparametric estimators are proposed. The estimators adapt to the regularity of the  process and to the nature of the design (sparse versus dense,  independent versus common).  They are simple  `first smooth, then estimate' procedures where the kernel estimates of the sample paths are constructed with optimal plug-in local bandwidths. The  bandwidths realize the minima of explicit pointwise risk bounds for the mean and autocovariance functions estimators, respectively. We also prove the pointwise asymptotic normality of the mean  estimator, a result which permits the construction of honest confidence intervals for non-differentiable mean functions. The study could be extended to other types of dependence, aiming at proving uniform convergence for the mean and (auto)covariance functions, or permitting for informative design \citep[see][]{weaver2022functional}. Such extensions are left for future work.




\appendix


\section{Technical lemmas}\label{app:example}

\subsection{$\LL^p-m-$approximability}

Let us  introduce some additional notations~: 
$	(f\otimes g)(s,t) = f(s)g(t)$, $\forall s,t \in I$  and $\ell \in \ZZ$.
Meanwhile, the tensor product $\circ$ is defined   as
$	\left(X_n \circ Y_n\right)(g)=\left<Y_n,g\right>_\Hh X_n$, for all $ X_n, Y_n, g\in \Cc$.
Recall that $\Ll = \Ll(\Cc, \Cc)$ is the space of bounded linear operators on $\Cc(I)$ equipped with the sup-norm. The proofs of the following  lemmas are given in \citet{hassan2024supp}.

\begin{lemma}\label{lem:lp-m-approx-1}
	\!\!\!\!	Let $\{X_n\} $,  $\{Y_n\} $ be  $\LL^p_\Cc\!-m-$approximable sequences, 
	for some $p\geq 4$. Define: 
	\begin{itemize}
		\item
		$Z_n^{(1)} = A(X_n)$, where $A \in \Ll$ ; and 	$Z_n^{(6)} = X_n\otimes X_{n+\ell}$, where here $\{X_n\} $ is $\LL^p_\Cc-m-$approximable for some $p\geq 8$.
		\item
		$Z_n^{(2)} = X_n + Y_n$; $Z_n^{(3)} = X_n  Y_n$; $Z_n^{(4)} = \langle X_n, Y_n \rangle_\Hh\in\RR$; and
		$Z_n^{(5)} = X_n \circ Y_n\in\Ll$.
	\end{itemize}
	Then $\{Z_n^{(1)}\}$, $\{Z_n^{(2)}\}$ are $\LL^p_\Cc-m-$approximable 
	in $\Cc$, 
	$\{Z_n^{(6)}\}$ is $\LL^{p/2}_\Cc-m-$approximable 
	in $\Cc$ and its  $\LL^{p/2}_\Cc-m-$approximation is $X_n^{(m)}\otimes X_{n + \ell}^{(m+\ell)}$. 
	If $X_n$ and $Y_n$ are independent, then $\{Z_n^{(3)}\}$, $\{Z_n^{(4)}\}$ and $\{Z_n^{(5)}\}$ are $\LL^p-m-$approximable in their respective spaces.
\end{lemma}

\begin{lemma}\label{lem:lp-m-approx-2}
	Let $\{X_n\}$ be a $\LL^p_\Cc-m$-approximable sequence. Let $s,t\in I$, $t \neq s$, and let $c$ be a constant. Define  $ F_n  = X_n(t)$ and
	$G_n  = \left(X_n(s) - X_n(t)\right)^{2}+c.$
	Then $\{F_n \}\subset\RR$ is $\LL^p-m-$approximable in $\LL^p$ 
	and $\{G_n\}\subset\RR$ is $\LL^{p/2}-m-$approximable in $\LL^{p/2}$.
\end{lemma}

\subsection{Study of \texorpdfstring{$\widehat{\theta}(u,v) $}~} 

Below  $(H_t,L_t^2)$ is a short notation for $(H_{0,t},L_{0,t}^2)$. The local regularity estimators are built as estimators of the proxy quantities $\widetilde H_t $ and $\widetilde L^2_t$. The concentration of our estimators will thus depend, on the one hand, on the accuracy  of the proxies, and on the other hand, on the concentration of the proxies' estimators based on the quantities  $\widehat\theta(u,v)$ defined in \eqref{def_theta_hat}. We  first investigate  these aspects before proving  Theorems~\ref{thm:regularity:H} and \ref{thm:regularity:L}.

\begin{lemma}[Proxies accuracy]\label{lem:convergence_proxy}
Let $t\in J$. 
\begin{enumerate}
	\item For any $\varphi\in(0,1)$ and $0<\Delta\leq \Delta_{0,0}$ such that $ 4\Delta^{2\beta_0} S_0^2 < L_t^2\log(2) \varphi$, we have 
	$
	| \widetilde H_t  - H_t|< \varphi/2.
	$
	\item  Let $H\in(0,1]$ such that $|H-H_t|<\varphi<1$.  
	For any $\psi\in(0,1)$ and $0<\Delta\leq \Delta_{0,0}$
	such that $S_0^2\Delta^{2\beta_0-2\varphi} < \psi/3$, we have
	$\Delta^{-2H}| \theta(t_1,t_3) - L_t \Delta^{2H_t}| < \psi/3$.
\end{enumerate} 
\end{lemma}

To study the properties of $\widehat\theta(u,v)$,   we use the Nagaev-type inequality   for sums of dependent random variables, see \citet{liu2013probability}. When dealing with real-valued variables, the dependence measure used under the $\LL^p-m-$approximability assumption is slightly more restrictive that the \textit{functional dependence measure}, as defined in \citet[Definition 1]{wu2005nonlinear}.  Lemma~\ref{lem:nagaev} below, a version of \citet[Theorem 2]{liu2013probability},  states a Nagaev-type inequality as we will use for our study. The proof is provided in \citet{hassan2024supp}. Let $\{U_n\}$ be a sequence of real-valued random variables, and let
$S_N^*=\max\{|S_n|, \ n = 1, \dots, N\}$,   where  $S_n = U_1 + \cdots + U_n$.

\begin{lemma}[Nagaev inequality]
\label{lem:nagaev}
Let $\{U_{n}\}\subset \RR$ be stationary,   $\LL^p-m-$approximable, 
\begin{equation*}
	\EE(U_{n})=0, \quad \text{and} \ \, \upsilon \!:=\! \sum_{m=1}^{\infty}\! \left( m^{p/2 - 1} \nu_{m,p}^p \right)^{1/(p+1)} \!\! < \infty \quad \text{ where}  \ \, 	\nu_{m,p} = \nu_p\left(U_m - U_{m}^{(m)}\right).
\end{equation*}
Then
\begin{equation*} 
	\PP\left( S_N^* \ge \varepsilon \right) 
	\le
	c_p\frac{N}{\varepsilon^p} \left( \upsilon^{p+1} + \|U_1\|_{p}^{p} \right) +
	c_p^{\prime}\exp{ \left( - \frac{c_p \varepsilon^2}{N \upsilon^{2+2/p} } \right) }
	+ 2\exp{ \left( - \frac{c_p \varepsilon^2}{N \|U_1\|_{2}^{2} } \right) },
\end{equation*}
where $c_p = 29p/\log(p)$ and $c_p^{\prime}$ are two positives constants.
\end{lemma}

We now study the concentration of $\widehat\theta(u,v)$ and $\widehat\theta(u,v)/\theta(u,v)$.

\begin{lemma}\label{lem:theta:conv}
Assume  the conditions  of Theorem~\ref{thm:regularity:H} hold true.
Let $u,v\in J$, $u\leq t\leq v$, be  such that $\Delta/2 \leq |u-v|\leq \Delta$ and let 
$\eta_0= \eta_0 (\lambda)= {8}\left(2\sqrt{\overline{a}_0} + \sqrt{R_2(\lambda)}\right)\sqrt{R_2(\lambda)}. $
For any $\kappa > 0$, define the probabilities
\begin{equation*}
	p_0^+(u,v; \kappa) = \PP\left[\widehat\theta(u,v)> (1+\kappa)\theta(u,v)\right],\quad 
	p_0^-(u,v;\kappa ) = \PP\left[\widehat\theta(u,v)< (1-\kappa)\theta(u,v)\right].
\end{equation*}
Then, for any $\eta$ such that  $\eta_0< \eta < 1$, 
\begin{equation*}
	\PP\left( \left| \widehat\theta(u,v) - \theta(u,v) \right| > \eta \right)
	\le
	\frac{\mathfrak{a}}{N\eta^2} +
	\mathfrak{b}\exp{ \left( - \mathfrak{e}N\eta^2 \right) },
\end{equation*}
where $\mathfrak{b}$ is a universal constant, and $\mathfrak{a}$ and $\mathfrak{e}$ are  constants  depending on the dependence measure and the fourth-order moment of $\widetilde X(u)$. 
Moreover, for any $\kappa$ such that $\eta_0<\kappa \theta(u,v)<1$, 
\begin{equation*}
	\max\big[p_0^+(u,v;\kappa ), p_0^-(u,v;\kappa )\big]
	\leq
	\frac{2^{2H_t+2}\mathfrak{a}}{N\kappa^2 {L_t^4} \Delta^{4H_t}} +
	\mathfrak{b}\exp{ \left( - \frac{\mathfrak{e}}{2^{H_t+2}}N\kappa^2 
		{L_t^4} \Delta^{4H_t} \right) }.
\end{equation*} 
\end{lemma}


\section{Proofs of Theorems~\ref{thm:regularity:H} and \ref{thm:regularity:L}}

\begin{proof}[Proof of Theorem \ref{thm:regularity:H}]
According to  \eqref{eq:thmH:proxy} and Lemma \ref{lem:convergence_proxy}, we have that $|\widetilde{H}_t-H_t|\leq \varphi/2.$ We then deduce
\begin{align*}
	\PP(|\widehat H_t - H_t| > \varphi)
	&\leq \PP\left(\left|\widehat H_t - \widetilde H_t \right| > \varphi/2\right)\leq \PP\left( \left|\log \frac{\widehat\theta(t_1,t_3)}{\theta(t_1,t_3)}\frac{\theta(t_1,t_2)}{\widehat\theta(t_1,t_2)}\right|  > \varphi\log(2)\right)\\
	&\leq \PP\left(\frac{\widehat\theta(t_1,t_3)}{\theta(t_1,t_3)}\frac{\theta(t_1,t_2)}{\widehat\theta(t_1,t_2)}  > 2^{-\varphi} \right) 
	+ \PP\left(\frac{\widehat\theta(t_1,t_3)}{\theta(t_1,t_3)}\frac{\theta(t_1,t_2)}{\widehat\theta(t_1,t_2)}  < 2^{-\varphi} \right).
\end{align*}
By simple algebra and using the definition of the functions $p_0^+$, $p_0^-$ from Lemma~\ref{lem:theta:conv}, we get
\begin{align*}
	\PP(|\widehat H_t - H_t| > \varphi)
	&\leq p_0^+(t_1,t_3; 2^{\varphi/2}-1) + p_0^-(t_1,t_3; 1-2^{-\varphi/2})\\
	&\qquad + p_0^+(t_1,t_2; 2^{\varphi/2}-1) + p_0^-(t_1,t_2; 1-2^{-\varphi/2}),
\end{align*}
provided that $\eta_0(\lambda) < |2^{\pm \varphi/2}-1| \theta(u,v)<1$. This is guaranteed by the condition \eqref{eq:thmH:tau-phi} with $C= B^{-1/2} (2\sqrt{\overline{a}_0} +\sqrt{B})^{-1}\log(2)/ 2^{15/2}$. To see this, first note that for any $\varphi \in (0,1)$, $|2^{\pm \varphi/2}-1| \leq \varphi \log(2) / 2^{1/2}$. By \eqref{eq:D:holder} and \eqref{eq:thmH:proxy}, we thus have
$$
|2^{\pm \varphi/2}-1| \theta(u,v) \leq \left(5 \log(2) / 2^{5/2}\right) \varphi L_t^2 \Delta^{2H_t} < 1 \quad \text{ as } \quad \Delta \rightarrow 0.
$$
Second, by \assrefH{H:pre-smooth:convergence},  $\eta_0 (\lambda) < {8}\left(2\sqrt{\overline{a}_0} + \sqrt{B}\right) B^{1/2}\lambda^{-\tau/2}$.
Gathering the two bounds, we get 
$$
\lambda^{-\tau/2} < \left(5 B^{-1/2} \left(2\sqrt{\overline{a}_0} + \sqrt{B}\right)^{-1} \log(2) / 2^{11/2}\right) \varphi L_t^2 \Delta^{2H_t},
$$
which is  condition \eqref{eq:thmH:tau-phi}.
Now, with $t_k = t_2$ or $t_k=t_3$, we have
\begin{multline*}
	p_0^\pm(t_1,t_k;\pm(2^{\pm \varphi/2}-1))
	\leq\!
	\frac{2^{2H_t+2}\mathfrak{a}}{N(2^{\pm\varphi/2}-1)^2 L_{t}^4 \Delta^{4H_t}} +
	\mathfrak{b}\exp\!\!\left(\! -\frac{\mathfrak{e}(2^{\varphi/2}-1)^2}{2^{2H_t+2}}N 
	L_{t}^4 \Delta^{4H_t}\! \right)\\
	\leq
	\frac{ 2^{2H_0+4}\mathfrak{a} / \log(2)^2}{N \varphi^2 L_t^4\Delta^{4 H_0} }
	+ \mathfrak{b} \exp{ \left( - 
		\frac{ \mathfrak{e} \log(2)^2 }{ 2^{2H_0+4} }  
		N \varphi^2  L_{t}^4 \Delta^{4 H_0} \right) }\\
	=	\frac{ \mathfrak{f}_0/4 }{ N \varphi^2 \Delta^{4 H_0} }
	+ \mathfrak{b} \exp{ \left( - \mathfrak{g}_0 N \varphi^2 \Delta^{4 H_0} \right) }	,
\end{multline*}
where $\mathfrak{f}_0 = 2^{2H_t+6}\mathfrak{a} / (\log(2)^2L_{t}^4)$ and $\mathfrak{g}_0 = \mathfrak{e}L_{t}^4\log(2)^2 / 2^{2H_t+4}$.
For the second inequality,  use  $ \log^2(2 )\varphi^2/4=\{\pm\log( 2^{\pm \varphi/2} ) \}^2\leq \{\pm (2^{\pm\varphi/2} - 1)\}^2$. 
The quantity $p_0^\mp(t_1,t_k;\mp(2^{\pm \varphi/2}-1))$ can be bounded similarly using   $ \log^2(2 )\varphi^2/4\leq \{\mp (2^{\pm\varphi/2} - 1)\}^2$. Gathering the four terms and changing $4\mathfrak b$ to $\mathfrak b$, we get the result. 
\end{proof}

\begin{proof}[Proof of Theorem \ref{thm:regularity:L}]
By definition and elementary algebra, 
$$
|\widehat L^2_t - {L}_t^2| \leq \frac{|\widehat{\theta}(t_1,t_3)-\theta(t_1,t_3)|}{\Delta^{2\widehat{H}_t}} + \frac{|\theta(t_1,t_3)-L_t^2\Delta^{2{H}_t}|}{\Delta^{2\widehat{H}_t}} + L_t^2|1-\Delta^{2H_t-2\widehat{H}_t}|,
$$
and thus  
\begin{align*}
	\PP(|\widehat L^2_t - L_t^2| > \psi)
	&\leq \PP\left(\left|\widehat H_t - H_t\right| \leq \varphi, |\widehat L^2_t - {L}_t^2| > \psi\right) + \PP\left(\left|\widehat H_t - H_t\right| > \varphi\right).
\end{align*}
If $\left|\widehat H_t - H_t\right| \leq \varphi$, by condition \eqref{eq:thmL:proxy} and Lemma~\ref{lem:convergence_proxy}, we get 
$\Delta^{-2\widehat{H}_t}|\theta(t_1,t_3)-L_t^2\Delta^{2{H}_t}| \leq \psi/3.$
Furthermore, since the function $x\to \Delta^{2x}$ is Lipschitz continuous over $[-\varphi,\varphi]$, we get
\begin{align*}
	L_t^2|1-\Delta^{2H_t-2\widehat{H}_t}|\leq \psi/3,
\end{align*}
provided $|\widehat H_t - H_t| \leq \varphi$ and condition  \eqref{eq:thmL:phi-psi} holds true.
We deduce that 
\begin{multline}
	\hspace{-.8cm}\PP(|\widehat L^2_t - L_t^2| > \psi) \leq 	\PP\left(|\widehat H_t - H_t| \leq \varphi, \left|\widehat{\theta}(t_1,t_3)-\theta(t_1,t_3)\right| > \Delta^{2\widehat{H}_t}\psi/3\right) \\\hspace{.4cm} + \PP\!\left(|\widehat H_t - H_t| > \varphi\right)
	\!\leq\! \PP\!\left(|\widehat{\theta}(t_1,t_3)-\theta(t_1,t_3)| > \Delta^{2H_t+2\varphi} \psi/3\right)\! + \PP\!\left(|\widehat H_t - H_t| > \varphi\right).
\end{multline}	
%
%
%
The second probability of the right-hand side of the last inequality can be bounded using Theorem \ref{thm:regularity:H} and the probability $p_0^+$ in Lemma \ref{lem:theta:conv}, provided that $\eta_0(\lambda) < \Delta^{2H_t+2\varphi} \psi/3 <1$ which is guaranteed by condition \eqref{eq:thmL:tau-psi} with $\widetilde{C}= B^{-1/2} (2\sqrt{\overline{a}_0} +\sqrt{B})^{-1} / (3 \times 2^{3})$. In fact, note that Assumption \assrefH{H:pre-smooth:convergence} implies $\eta_0 (\lambda) \leq {8}\left(2\sqrt{\overline{a}_0} + \sqrt{B}\right) B^{1/2}\lambda^{-\tau/2}$, hence
$$
\lambda^{-\tau/2} < B^{-1/2} (2\sqrt{\overline{a}_0} +\sqrt{B})^{-1}/ (3 \times 2^{3}) \Delta^{2\varphi} \psi \Delta^{2H_t} < 1\quad \text{ as } \quad \Delta \rightarrow 0.
$$
Then Theorem \ref{thm:regularity:L} follows. 
\end{proof}


\section{Adaptive estimation} \label{sec_adapt_app}

\subsection{Technical lemmas}
Let $X$ be a generic random function having the stationary distribution
of $\{X_n\}$. Let  $\EEn(\cdot)\! = \!\EE(\cdot \mid \! M_n,\!\{\Tni, 1\leq i \leq M_n\!\}, X_n\!)$ and $\EEMT(\cdot) = \EE \left(\cdot \mid \! M_n,\! \{\Tni, 1\leq i \leq M_n\!\}, 1\leq n \leq N \right)$.
Below, `wrt' is the abbreviation for `with respect to'. We consider the  decomposition
$\widehat X_n(t;h) - X_n(t) = B_n(t;h) + V_n(t;h) $, $t\in I$,
into a bias term and a stochastic term, where 
$$
B_n(t;h) \coloneqq \EEn\left[ \widehat X_n(t;h) \right]  -  X_n(t) 
\quad\text{ and  }\quad
V_n(t;h) \coloneqq \widehat X_n(t;h) - \EEn \left[ \widehat X_n(t;h) \right]. 
$$
By construction due to \assrefH{H:epsni} and \assrefH{H:ind}, $\forall n\neq n^\prime$, $\EEMT\![V_n(t;h)]\!=\!0,$ $\EEMT\![V_n(t;h)V_{n^\prime}(t;h)]\!=\!0,$ $\EEMT[V_n(t;h)B_n(t;h)]=0$ and $\EEMT[V_n(t;h)B_{n^\prime}(t;h)]  = 0$. 
The following resultstudies  $\{B_n(t;h)\}$ and $\{V_n(t;h)\}$. Here, $\{B_n(t;h)\}$ is a short notation for $\{B_n(t;h), 1\leq n\leq N\}$, and the same rule is used for $\{V_n(t;h)\}$ and $\{M_n\}$, while $\{\Tni\}$ means $\{\Tni, 1\leq n \leq N, 1\leq i\leq M_n\}$. Below, $\Xx( H , \boldsymbol L ;J)$ is the class from Definition \ref{def:space_family}, with $\delta=0$. 

\begin{lemma} \label{lem:biais-var} 
Assume that $X\in \Xx( H , \boldsymbol L ;J)$ and let $t\in I$ and  $\widehat X_n(t,h)$ be defined as in \eqref{LP_est_v}. 
Assume \assrefH{H:stationarity} to \assrefH{D:smooth2}, \assrefH{H:polynomelocaux} and \assrefH{H:cd_grid_H}.
Then~:   
\begin{enumerate}
	\item\label{mean_l1}  
	$\{B_n(t;h)\}$ and $\{V_n(t;h) \}$ are conditionally independent given  $\{M_n\}$ and $\{\Tni\}$~;
	
	\item\label{mean_l2} 
	$\{V_n(t;h) \}$ are conditionally independent given $\{M_n\}$ and $\{\Tni\}$~;
	
	\item\label{mean_l3} 
	$\EEMT \!\left[V^2_{n}(t;h)\right]\!  \leq \! \{1+o(1)\} \sigma^2 (t)  \max_{1\leq i \leq M_n}\!W_{n,i}(t;h)$, with $o(1)$ uniform wrt $h\in\mathcal H_N$;
	
	\item\label{mean_l4} 
	$\EEMT \left[B^2_{n}(t;h)\right] \leq L_t^2 h^{2H_t}  b_n(t;h,2H_t)\{1+o(1)\}$, with $o (1)$  uniform wrt $h\in\mathcal H_N$.

\end{enumerate}
\end{lemma}

\begin{lemma}\label{lem:concen-PN-PNl}
Assume that the assumptions \assrefH{H:stationarity} to  \assrefH{H:ind}, and \assrefH{H:cd_grid_H} to \assrefH{H:densityT} hold true.
\begin{enumerate}
	\item\label{lemP_1} 
	For all $t\in(0,1),$ and $h\in\Hh_N$  
	\begin{equation*}
		1- \exp\left(-M p(t;h) \right)  \leq \EE[\pi(t;h)|M_1] \leq 1-\exp\left(-2M  p(t;h)\right)\quad \text{a.s}.
	\end{equation*}
	
	\item\label{lemP_2} 
	There exists two constants $\underline{C}_\mu$ and $\overline{C}_{\mu}$ such that for all $h\in\Hh_N$,
	\begin{equation*}
		\underline{C}_\mu \{1 + o(1)\} \leq \frac{\EE[P_N(t;h)]}{N\min(1,\lambda h)}
		\leq  \overline{C}_{\mu} \{1+o(1)\},
	\end{equation*}
	and $P_N(t;h) = \EE[P_N(t;h)]\{1 + o_\PP(1)\},$ with $o(1)$ and $o_\PP(1)$ uniform wrt  $h\in \Hh_N$.
	
	\item\label{lemP_3} 
	Moreover if \assrefH{H:cd_grid_H_autocov} holds, constants $\underline{C}_{\gamma}$ and $\overline{C}_{\gamma}$ exist such that  $\forall h\in\Hh_N$,
	$$
	\underline{C}_{\gamma} \{1 + o(1)\}\leq \frac{\EE[P_{N,\ell}(s,t;h)]}{(N-\ell)\min(1, (\lambda h)^2)} 
	\leq  \overline{C}_{\gamma} \{1+o(1)\},
	$$
	and $P_{N,\ell}(s,t;h) = \EE[P_{N,\ell}(s,t;h)]\{1 + o_\PP(1)\},$ with $o(1)$ and $o_\PP(1)$ uniform  wrt  $h\in \Hh_N$.
\end{enumerate}
\end{lemma}

\begin{lemma}\label{lem:cv-sigma}
\!\!If assumptions \assrefH{H:stationarity} to \assrefH{H:Lpmapprox} and \assrefH{H:densityT} hold true,
$\widehat\sigma^2(t)= \sigma^2(t)\{1 + o_\PP(1)\}$.  
\end{lemma}

\begin{lemma}\label{NWweigths-bound} 
Assume the assumptions \assrefH{H:stationarity} to \assrefH{H:ind}, and \assrefH{H:polynomelocaux} to \assrefH{H:densityT} hold true.
For each 	$N\geq 1$, we have
$$
0\leq \max_{n,i} W_{n,i}(t;h) \leq S_{n,W}(h)\min\left(1, (\lambda h)^{-1}\right), \quad 1\leq n \leq N,
$$
where $S_{n,W}(h)\geq 1$ is a random variable with the mean and the variance bounded by constants.  Moreover, the variables  $\{S_{n,W}(h), 1 \leq n\leq N\}$ are independent.  
\end{lemma}

\begin{lemma}\label{lem:cv-wsum}
Assume that the assumptions \assrefH{H:stationarity} to \assrefH{H:ind}, \assrefH{H:Lpmapprox} for $p \geq 8$ and \assrefH{H:cd_grid_H} hold. 
For each $h\in\Hh_N$,  let $\{\pi_n(h), n\geq 1\}$ be a sequence of i.i.d. Bernoulli random variables which is independent of $\{X_n, n\in \ZZ\}$.
Then, for any $t \in I $, 
$N^{-1} \sum_{n=1}^{N} \pi_n(h)X^2_n(t) = \EE\left[\pi_1(h)X^2_1(t)\right]\{1 + o_\PP(1)\}$ uniformly wrt  $h \in \Hh_N$.
\end{lemma}

\subsection{\!Mean function: risk bound, rates of convergence, asymptotic normality}\label{sec:mean_ad_rate}

\begin{lemma}\label{lem:mu:risk-bound} 
Under the assumptions \assrefH{H:stationarity} to \assrefH{H:Lpmapprox}, \assrefH{H:polynomelocaux} and \assrefH{H:cd_grid_H}, we have
\begin{equation*}
	\EEMT\left[\{\widehat \mu_N(t;h) - \mu(t)\}^2\right] \leq 2R_\mu(t;h)   \{1+ o (1)\}, 
\end{equation*}
with $o (1)$  uniform wrt $h\in\mathcal H_N$ and  
\begin{equation*}
	R_\mu(t;h) = L_t^2 h^{2 H_t}\mathbb{B} (t;h,2H_t) + \sigma^2(t) \mathbb{V}_\mu(t;h) + \mathbb{D}_\mu (t;h)/P_N(t;h).
\end{equation*}

\end{lemma}

\begin{lemma}\label{lem:risk-bound:cv} 
Assume the assumptions \assrefH{H:stationarity} to \assrefH{H:Lpmapprox}, \assrefH{H:cd_grid_H}, \assrefH{H:densityT}, \assrefH{H:HL_conv} hold true. 
Let 
\begin{equation*}
	\widehat R_\mu(t;h) = \widehat L_t^2 h^{2\widehat H_t}\mathbb{B} (t;h,2\widehat H_t) +  \widehat\sigma^2(t)   \mathbb{V}_\mu(t;h) + \mathbb{D}_\mu(t;h)/P_N(t;h).
\end{equation*}
Then 
$
\sup_{h\in\Hh_N}  \widehat R_\mu(t;h)/R_\mu(t;h)   = 1+ o_\PP (1).
$
\end{lemma}

\begin{proof}[Proof of Theorem~\ref{thm:mu}: rates of convergence]
We recall that 
\begin{equation*} 
	\widehat R_\mu(t;h) =   \widehat L_t^2 h^{2 \widehat H_t}\mathbb{B} (t;h,2 \widehat H_t) +   \widehat  \sigma^2(t)  \mathbb{V}_\mu(t;h) + {\mathbb{D}}_\mu(t;h)/P_N(t;h).
\end{equation*}
Let us define  $\mathcal H_{1,N} = \{h\in \mathcal H_N, \lambda h \leq C\}$ and $\mathcal H_{2,N} = \{h\in \mathcal H_N, \lambda h > C\}$, for some $C\geq 1$. 
Thus, we have  $\mathcal H_{N}=\mathcal H_{1,N}\cup \mathcal H_{2,N}$. 
Over the set $\mathcal H_{1,N} $, we simply recall  $\displaystyle \max_{1\leq i \leq M_n}\!  W_{n,i}(t;h) \leq 1$, for any $h\in\mathcal H_{1,N}$. 
Over the set $\mathcal H_{2,N} $, 
Lemma~\ref{NWweigths-bound} implies
\begin{align*}
	&\mathbb{V}_\mu(t;h) \leq   \frac{ \min\left(1, (\lambda h)^{-1}\right)}{P _N(t;h) } \times \dfrac{1}{P_N(t;h)} \sum_{n=1}^N  \pi_n(t;h)S_{n,W}(h) ,\\
	\text{with }& S_{n,W}(h) =\max \left\{1  ,\; (\lambda h) (\|K\|_\infty/\tau)s_{n,W}(h)  \right\},\quad s_{n,W}(h) =  \frac{\mathds{1}\{S(M_n,t,h)>0 \}}{S(M_n,t,h) },	
\end{align*}
and $S(M_n,t,h)$  the integer-valued variable,  non-decreasing as function of $h$, which counts the number of points in $[t-h,t+h]$. 
Since, for  $a\geq 0$ we have $\max\{1,a\}\leq 1+a$, we  study 
$$
\sum_{n=1}^N  \pi_n(t;h)s_{n,W}(h) = \sum_{n=1}^N  s_{n,W}(h).
$$
The equality is the consequence of the fact that, by definition,  
$
\mathds{1}\{S(M_n,t,h)>0 \} \pi_n(t;h) = \mathds{1}\{S(M_n,t,h)>0 \}. 
$
By the calculations provided in the proof of Lemma~\ref{NWweigths-bound}, 
$$
c_{1,W} / (\lambda h) \leq 	\EE \left[   s_{n,W}(h) \right]\leq c_{2,W} / (\lambda h)\quad \text{and} \quad  \operatorname{Var}( s_{n,W}(h) )\leq \EE \left[  s^2_{n,W}(h) \right] \leq c^\prime _W/ (\lambda h)^2,
$$
where $c_{1,W},c_{2,W}$ and $c^\prime_W$ are positive constants depending only on $C$, $\underline c_g$,  $\overline c_g$, $\underline c$ and $\overline c$. Moreover, $\{s_{n,W}(h), n\geq 1\}$ is a sequence of independent variables, bounded by 1. Applying Bernstein's inequality \citep[see][Theorem 1.8.4]{V_2018} for each $h\in\mathcal H_{N,2}$,  we deduce that $\forall \epsilon >0$, a constant 
$C_\epsilon >0$ exists, depending on $\epsilon$,  $C$,
$c_{2,W}$ and $c^\prime_W$, such that 
$$
\PP \left(	\dfrac{1}{N} \sum_{n=1}^N    s_{n,W}(h) > \EE \left[   s_{n,W}(h) \right] + \epsilon  \right)  \leq \exp\left( - C_\epsilon N \right)  .
$$
Using next 	Boole's (union bound) inequality, at the price of a logarithmic term in the exponential, we get a uniform over $\mathcal H_{2,N}$ exponential  bound  for the upper tail probability of $N^{-1}\sum_{n=1}^N s_{n,W}(h)$. By Lemma \ref{lem:concen-PN-PNl}-(\ref{lemP_2})  we deduce that  a constant $\mathfrak c>0$ exists such that 
$$
\sup_{h\in \mathcal H_{2,N}}\dfrac{1}{P_N(t;h)} \sum_{n=1}^N  \pi_n(t;h)S_{n,W}(h)\leq  \mathfrak c+o_\PP (1),
$$ 
and thus
\begin{equation}\label{rez_V}
	\mathbb{V}_\mu(t;h) \leq \min\left(1, (\lambda h)^{-1}\right) \times \{P _N(t;h) \}^{-1} \times\{\mathfrak c+o_\PP (1)\},
\end{equation} 
uniformly over $\mathcal H_{N}$.  According to \assrefH{H:HL_conv},  $\widehat L_t^2 = L_t^2 \{1 + o_\PP(1)\}$.  
Moreover,  by Lemma \ref{lem:risk-bound:cv}, $h^{2\widehat H_t} = h^{2H_t} \{1 + o_\PP(1)\}$ uniformly over the grid $\Hh_N$, and, by  Lemma~\ref{lem:cv-sigma},  $\widehat{\sigma}^2(t)=\sigma^2(t) \{1 + o_\PP(1)\}.$
From these, and simply bounding $\mathbb{B} (t;h,\alpha)$ by 1, we get 
\begin{equation*}
	\widehat R_\mu(t;h) \leq L_t^2 h^{2 H_t} \{1 + o_\PP(1)\} + \sigma^2(t) \{\mathfrak c + o_\PP(1)\} \dfrac{\min\left(1,\; (\lambda h)^{-1}\right)}{P_N(t;h)} + \dfrac{\mathbb{D}_\mu(t;h)}{P_N(t;h)},
\end{equation*}
uniformly over $\mathcal H_{N}$.
By Lemma~\ref{lem:concen-PN-PNl}, we have
$$
\dfrac{1}{P_N(t;h)}\leq \dfrac{\underline{C}_\mu^{-1}}{N\min(1, \lambda h)}\{1 + o_\PP(1)\},
$$
and 
\begin{equation*}
	\dfrac{\min\left(1,\; (\lambda h)^{-1}\right)}{P_N(t;h)}
	\leq \underline{C}_\mu^{-1}\dfrac{\min\left(1,\; (\lambda h)^{-1}\right)}{N\min(1, \lambda h)}\{1 + o_\PP(1)\} = \dfrac{\underline{C}_\mu^{-1}}{N \lambda h}\{1 + o_\PP(1)\} ,
\end{equation*}
with the $o_\PP(1)$ uniform with respect to $h\in\Hh_N$. Gathering facts, we get 
\begin{equation*}
	\widehat R_\mu(t;h) =  \Oo_\PP\left\{h^{2H_t} + (N \lambda h)^{-1} + N^{-1}\right\}.
\end{equation*}
The right-hand side is minimized by $h$ with the rate $(N \lambda )^{-1 / \{2H_t + 1\}}$. The rate for convergence of $\widehat\mu_N^*(t) - \mu(t)$ is obtaining by replacing the optimal bandwidth in the risk bound. Indeed, in the case where $\lambda ^{2H_t} \ll N$ (sparse case), we get 
$
(N \lambda )^{-2H_t / \{2H_t + 1\}} \gg N^{-1}, 
$
and $\widehat\mu_N^*(t) - \mu(t)$ converges at the rate $O_\PP((N \lambda )^{-H_t / \{2H_t + 1\}})$, which is slower than $O_\PP( N^{-1/2})$. Meanwhile,
when $\lambda ^{2H_t} \gg N$ (dense case), we have 
$
(N \lambda )^{-2H_t / \{2H_t + 1\}} \ll N^{-1}, 
$
and the rate of convergence of $\widehat\mu_N^*(t) - \mu(t)$ is given by the square root of 
$\mathbb{D}_\mu(t;h)/P_N(t;h)$ which in this case leads to the parametric rate $O_\PP( N^{-1/2})$. \end{proof}

\begin{proof}[Proof of the Theorem~\ref{thm:mean-clt}]
Let 
$\widetilde \mu_N (t;h) = \{P_N(t;h) \}^{-1}\sum_{n=1}^N  \pi_n (t;h) X_n(t) $.
Then,
\begin{equation*}
	\widehat \mu_N(t;h) - \mu(t) = \left\{\widehat \mu_N(t;h) - \widetilde  \mu_N(t;h)\right\}+ \left\{\widetilde \mu_N(t;h) -   \mu(t)\right\}  =:G_{N1}(t;h) + G_{N2}(t;h). 
\end{equation*}

\textbf{Convergence of $G_{N2}(t;h)$.} Given the indicators $\pi_n (t;h)$, $n\geq 1$, we use the CTL given in \citet[Theorem 3]{wu2011asymptotic} under predictive dependence. That result can be applied because, on the one hand, \citet[Theorem 1]{wu2005nonlinear} states that the functional (also called physical dependence) implies the predictive dependence, and, on the other hand, \citet[Lemma 1]{chen2015simultaneous} show that the $\LL^p-m$-approximation implies the functional dependence.
Then 
conditionally on 
$
\mathcal{X}_N = \{M_n, T_{n,i}, 1 \leq i \leq M_n, 1 \leq n \leq N\}
$
such that $P_N(t;h)\rightarrow \infty$ and 
\begin{equation}
	\mathbb{S}_{N,\mu} (t) :=\operatorname{Var}_{M,T} \left( \frac{1}{\sqrt{P_N(t;h)} }\sum_{n=1}^N  \pi_n (t;h) \{X_n(t)- \mu(t)\} \right),
\end{equation}
has a  limit in $(0,\infty)$, we have
\begin{equation}\label{tcl_1}
	\sqrt{P_N(t;h)/\mathbb{S}_{N,\mu} (t)}	\ G_{N2}(t) \overset{d}{\longrightarrow} \mathcal{N}\left(0, 1   \right).
\end{equation}

\textbf{Convergence of $G_{N1}(t;h)$.} Using \eqref{eq:data-model} and \eqref{LP_est_v}, we consider the decomposition
\begin{equation}
	G_{N1}(t;h) = \sum_{n=1}^N \frac{\pi_n (t;h) B_n(t;h)}{P_N(t;h)} + \sum_{n=1}^N  \frac{\pi_n (t;h) V_n(t;h)}{P_N(t;h)} := \mathcal{B}_N(t;h) + \mathcal{V}_N(t;h),
\end{equation}
where 
\begin{equation}
	B_n(t;h) = \sum_{i=1}^{M_n}W_{n,i}(t;h)\{X_n(T_{n,i}) - X_n(t)\} 
	\quad \text{and} \quad
	V_n(t;h) = \sum_{i=1}^{M_n}W_{n,i}(t;h)\sigma(T_{n,i})\varepsilon_{n,i}.
\end{equation}
We learn from the proof of Theorem~\ref{thm:mu} that any bandwidth sequence $h_N$ with a faster decrease than $(N \lambda)^{-1/(2H_t+1)}$ makes $\mathcal{B}_N $  negligible 
with respect to $\mathcal{V}_N$. This happens under the condition $h_N   (N\lambda)^{1/(2H_t+1)} \rightarrow 0$. We thus only have to study 
$\mathcal{V}_N(t;h)$. We can write 
\begin{equation}
	\mathcal{V}_N(t;h) 
	=  \frac{\{1+ o(h)\}\sigma(t)}{P_N(t;h)}  \sum_{n=1}^N  \pi_n (t;h) \sum_{i=1}^{M_n}W_{n,i}(t;h)\varepsilon_{n,i} =: \frac{1+ o(h)}{\sqrt{P_N(t;h)}}  \times \mathcal U_N (t;h).
\end{equation}
By Lyapunov CLT for independent variables, conditionally given the $M_n$ and $\{T_{n,i},  1 \leq i \leq M_n\}$, $1\leq n \leq N$, we have $A_N (t;h_N)^{-1/2}\mathcal U_N (h_N)\overset{d}{\longrightarrow} N(0,1)$ with 
\begin{equation*}
	A _N (t;h_N)=  \frac{\sigma^2(t)}{P_N(t;h_N)} \sum_{n=1}^N  \pi_n (t;h_N) \sum_{i=1}^{M_n}W^2_{n,i}(t;h_N).
\end{equation*}
This implies that for any sequence $A _N (t;h_N)$ which convergences to $\Sigma(t)$, we get 
\begin{equation}\label{zelm}
	\EE_{M,T}\left[\exp\left\{-iu\sqrt{P_N(t;h_N)}\;\mathcal{V}_N(t;h_N)\right\}\right] 
	\longrightarrow \exp(-u^2\Sigma(t) / 2) ,\qquad \forall u\in\RR.
\end{equation}
Note that $\EEMT[\cdots]$ on the left hand side is a bounded sequence of random variables. Since $A _N (t;h_N) - \Sigma(t) = o_\PP(1)$, and the convergence in probability is characterized by the fact that 
every sub-sequence  has a further sub-sequence which convergences almost surely, we deduce that the convergence in \eqref{zelm} holds in probability. By the Dominated Convergence Theorem for a sequence of bounded random variables 
convergent in probability, we get
\begin{equation}\label{zelm2}
	\EE \left[\exp\left\{-iu\sqrt{P_N(t;h_N)}\; \mathcal{V}_N(t;h_N) \right\}\right] 
	\longrightarrow \exp(-u^2\Sigma(t) / 2) ,\qquad \forall u\in\RR,
\end{equation}
which means 
\begin{equation}\label{tcl_2}
	\sqrt{P_N(t;h_N)}\; \mathcal{V}_N(t;h_N) \overset{d}{\longrightarrow} \mathcal{N}\left(0, \Sigma(t) \right).
\end{equation}
By   \assrefH{H:ind},  $G_{N2}(t)$ and $\mathcal{V}_N (t;h_N)$ are independent. From this, \eqref{tcl_1} and \eqref{tcl_2}, we get  
$$
\sqrt{P_N(t;h_N)}\; \{\mathcal{V}_N(t;h_N) + G_{N2}(t) \} \overset{d}{\longrightarrow} \mathcal{N}\left(0, \Sigma(t)+\mathbb{S}_\mu (t)  \right).
$$
By Lemmas \ref{lem:biais-var}-(\ref{mean_l3}) and \ref{NWweigths-bound}, $\Sigma(t)=0$ if $\lambda h_N \rightarrow \infty$, and $\Sigma(t)=\sigma^2(t)$ if  $\lambda h_N \rightarrow 0$.

Let us note that,
\begin{equation*}
	\mathbb{S}_{N,\mu} (t)
	=  \EE\left[\{X_0(t) - \mu(t)\}^2\right] 
	+ 2 \sum_{\ell=1}^{N-1}\EE\left[ \{X_0(t) - \mu(t)\}  \{X_\ell(t) - \mu(t)\} \right] \frac{P_{N,\ell}(t,t;h)}{P_N(t;h)},
\end{equation*}
with $P_{N,\ell}(s,t;h)$ defined in \eqref{PNl_def}.
It is easy to show that the conditional distribution of $\pi_n(t;h)\pi_{n+\ell}(t;h)$ given $P_N(t;h)$ is the Bernoulli distribution of success probability parameter  $P_N(t;h)(P_N(t;h)-1)\{N(N-1)\}^{-1}.$ We then get
\begin{align*}
	\EE\left[\mathbb{S}_{N,\mu} (t) \mid P_N(t;h)\right] =&  \EE\left[\{X_0(t) - \mu(t)\}^2\right]\\ & 
	+ 2 \frac{P_N(t;h)-1}{N-1} \sum_{\ell=1}^{N-1}\EE\left[ \{X_0(t) - \mu(t)\}  \{X_\ell(t) - \mu(t)\} \right] \frac{N-\ell}{N}.
\end{align*}
Using the Dominated Convergence Theorem, we get
$
\EE\left[\mathbb{S}_{N,\mu} (t) \mid P_N(t;h)\right] \rightarrow \mathbb{S}_{\mu} (t)$, in probability.
Lemma~\ref{lem:concen-PN-PNl} then implies that a constants $\underline C,\overline C$ exist such that 
\begin{multline}
	\!\!\!\!\!\!\!\!\!\!\!\!\!\!\EE\!\left[\{X_0(t) - \mu(t)\}^2\right] 
	+ 2\underline C\min(1,\lambda h) \sum_{\ell\geq 1}\EE\!\left[ \{X_0(t) - \mu(t)\}  \{X_\ell(t) - \mu(t)\} \right] \leq 	\mathbb{S}_{\mu} (t) \\\; \; \leq  \EE\!\left[\{X_0(t) - \mu(t)\}^2\right] 
	+ 2\overline C\min(1,\lambda h) \sum_{\ell\geq 1}\EE\!\left[ \{X_0(t) - \mu(t)\}  \{X_\ell(t) - \mu(t)\} \right].	 
\end{multline}
In particular, this means $\mathbb{S}_{\mu} (t) = \operatorname{Var}(X(t))$ provided $\lambda h_N \rightarrow 0$.  \end{proof}

\subsection{\!Autocovariance estimator: risk bounds and rates of convergence}

\begin{lemma}\label{lem:gamma:risk-bound} 
Under the assumptions \assrefH{H:stationarity} to \assrefH{D:smooth2}, \assrefH{H:Lpmapprox} for $p \geq 8$, \assrefH{H:polynomelocaux} to \assrefH{H:densityT}, and \assrefH{H:cd_grid_H_autocov} we have
$\EEMT\left[\left\{\widehat \gamma_{N,\ell}(s,t;h) -   \gamma_\ell (s,t)   \right\}^2 \right] \leq 2R_\gamma(s,t;h)  	\{ 1+o_\PP (1)\} $,
with $o_\PP (1)$  uniform with respect to $h\in\mathcal H_N$ and $	R_\gamma(s,t;h) $ defined in \eqref{R_ga_theo}.
	\end{lemma}

	\begin{proof}[Proof of Theorem~\ref{thm:gamma}]
By construction, $\mathbb{B} (t|s;h,\alpha,\ell^\prime )\leq 1$. Lemma \ref{NWweigths-bound} entails 
$$
\max \left\{ \mathbb{V}_{\gamma,0}(s,t;h), \mathbb{V}_{\gamma,\ell }(s,t;h)\right\}  \leq \{1+o_\PP(1)\}C_W\min\left(1,(\lambda h)^{-1}\right)/P_{N,\ell}(s,t;h),
$$
with $o_\PP(1)$ independent of $h$. See also the arguments used for \eqref{rez_V}.
By similar arguments
$$
\mathbb{V}_{\gamma}(s,t;h)  \leq \{1+o_\PP(1)\}
	C^2_W\min\left(1,(\lambda h)^{-2}\right)/P_{N,\ell}(s,t;h),
$$ 
uniformly with respect to $h$. 	
Next, from Lemma~\ref{lem:concen-PN-PNl}, we get
\begin{align}
	\dfrac{\min\left\{1, (\lambda h)^{-1}\right\}}{P_{N,\ell}(s,t;h)}& = 
	\dfrac{\min\left\{1, (\lambda h)^{-1}\right\}}{\min\{1, (\lambda h)^2\}}\times 
	\frac{\min\{ 1, (\lambda h)^2\} }{\EE[P_{N,\ell}(s,t;h)]} \times
	\frac{\EE[P_{N,\ell}(s,t;h)]}{P_{N,\ell}(s,t;h)}\nonumber\\
	& \leq \underline{C}_\gamma^{-1}\times [N\min\left\{\lambda h, (\lambda h)^2\right\}]^{-1}  \times
	\{1 + o_\PP(1)\} \label{eq:bound-ratio-PNl-1}.
\end{align}
By  similar calculations, 
\begin{equation}\label{eq:bound-ratio-PNl-2}
	\{P_{N,\ell}(s,t;h) \} ^{-1} \times  \min\left\{1, (\lambda h)^{-2}\right\} 
	\leq \underline{C}_\gamma^{-1}\times 
	[N\min\left\{1, (\lambda h)^2\right\}]^{-1}  
	\times\{1 + o_\PP(1)\}, 
\end{equation}
with $o_\PP(1)$ terms uniform with respect to $h\in\Hh_N$. Moreover, again using Lemma~\ref{lem:concen-PN-PNl}, we have 
\begin{equation}\label{eq:bound-PNl-3}
	\{P_{N,\ell}(s,t;h) \} ^{-1}
	\leq \underline{C}_\gamma^{-1} [N\min\left\{1, (\lambda h)^2\right\}]^{-1}  \times \{1 + o_\PP(1)\},
\end{equation}
uniformly with respect to $h\in\mathcal H_N$.

Let us now recall that from \assrefH{H:HL_conv} we have $\widehat L_t^2 = L_t^2 \{1 + o_\PP(1)\}$, uniformly over $h\in\mathcal H_N$. Moreover,  following   the lines of the proof of Lemma~\ref{lem:risk-bound:cv} we have $h^{2\widehat H_t} = h^{2H_t} \{1 + o_\PP(1)\}$ uniformly over the grid $\Hh_N$. Finally,  Lemma~\ref{lem:cv-sigma} establishes that $\widehat{\sigma}^2 (t)=\sigma^2(t) \{1 + o_\PP(1)\}$ uniformly. Therefore, $\widehat R_\gamma(s,t;h) =  R_\gamma(s,t;h) \{1 + o_\PP(1)\}$,  uniformly over $h\in\mathcal H_N$, with $R_\gamma(s,t;h)$ defined in \eqref{R_ga_theo}.
Gathering facts, and using equations \eqref{eq:bound-ratio-PNl-1}, \eqref{eq:bound-ratio-PNl-2} and \eqref{eq:bound-PNl-3}, we obtain
\begin{equation} \label{eq:bound-risk-rate-gamma}
	\widehat R_\gamma(s,t;h) = \Oo_\PP(h^{2H(s,t)} + \{N\min(\lambda h, (\lambda h)^2)\}^{-1}+  \{N\min(1, (\lambda h)^2)\}^{-1}  ),
\end{equation}
where $H(s,t) = \min(H_s, H_t)$.
The right-hand side of  \eqref{eq:bound-risk-rate-gamma} is minimized by a bandwidth
$	h_\gamma^* \sim \max\{ (N\lambda^2)^{-1/\{2H(s,t) + 2\}}, (N\lambda)^{-1/\{2H(s,t) + 2\}} \}$.
Finally, by replacing this rate in the equation~\eqref{eq:bound-risk-rate-gamma} we have the following rate of convergence
\begin{equation}
	\widehat \gamma_{N,\ell}(s,t;h_\gamma^*) - \gamma_\ell(s,t) 
	= \Oo_\PP\left((N\lambda^2)^{-\frac{H(s,t)}{2\{H(s,t) + 1\}}} + (N\lambda)^{-\frac{H(s,t)}{2H(s,t) + 1}} + N^{-1/2}  \right).
\end{equation}
\end{proof}



\bibliographystyle{apalike}
\bibliography{bibliography}

\end{document}